\documentclass[reqno,12pt]{amsart}
\usepackage{amscd,amssymb}

\theoremstyle{plain}
\newtheorem{Thm}[subsection]{Theorem}
\newtheorem{Cor}[subsection]{Corollary}
\newtheorem{Lem}[subsection]{Lemma}
\newtheorem{Prop}[subsection]{Proposition}
\newtheorem{Conj}[subsection]{Conjecture}

\theoremstyle{definition}
\newtheorem{Def}[subsection]{Definition}

\theoremstyle{remark}

\newtheorem{Rem}[subsection]{Remark}

\newtheorem{conj}{Conjecture}

\numberwithin{equation}{section}
\renewcommand{\rm}{\normalshape}

\newif\ifShowLabels
\ShowLabelstrue
\newdimen\theight
\def\TeXref#1{%
        \leavevmode\vadjust{\setbox0=\hbox{{\tt
                \quad\quad  {\small \rm #1}}}%
        \theight=\ht0
        \advance\theight by \lineskip
        \kern -\theight \vbox to
        \theight{\rightline{\rlap{\box0}}%
        \vss}%
        }}%

\ShowLabelsfalse

\renewcommand{\sec}[2]{\section{#2}\label{S:#1}%
        \ifShowLabels \TeXref{{S:#1}} \fi}
\newcommand{\ssec}[2]{\subsection{#2}\label{SS:#1}%
        \ifShowLabels \TeXref{{SS:#1}} \fi}

\newcommand{\sssec}[2]{\subsubsection{#2}\label{SSS:#1}%
        \ifShowLabels \TeXref{{SSS:#1}} \fi}

\newcommand{\refs}[1]{Section ~\ref{S:#1}}
\newcommand{\refss}[1]{Section ~\ref{SS:#1}}

\newcommand{\reft}[1]{Theorem ~\ref{T:#1}}
\newcommand{\refl}[1]{Lemma ~\ref{L:#1}}
\newcommand{\refp}[1]{Proposition ~\ref{P:#1}}
\newcommand{\refc}[1]{Corollary ~\ref{C:#1}}

\newcommand{\refe}[1]{\eqref{E:#1}}
\newcommand{\refco}[1]{Conjecture ~\ref{Co:#1}}

\newenvironment{thm}[1]%
        { \begin{Thm} \label{T:#1}  \ifShowLabels \TeXref{T:#1} \fi }%
        { \end{Thm} }

\renewcommand{\th}[1]{\begin{thm}{#1} \sl }
\renewcommand{\eth}{\end{thm} }

\newenvironment{lemma}[1]%
        { \begin{Lem} \label{L:#1}  \ifShowLabels \TeXref{L:#1} \fi }%
        { \end{Lem} }
\newcommand{\lem}[1]{\begin{lemma}{#1} \sl}
\newcommand{\elem}{\end{lemma}}

\newenvironment{propos}[1]%
        { \begin{Prop} \label{P:#1}  \ifShowLabels \TeXref{P:#1} \fi }%
        { \end{Prop} }
\newcommand{\prop}[1]{\begin{propos}{#1}\sl }
\newcommand{\eprop}{\end{propos}}

\newenvironment{corol}[1]%
        { \begin{Cor} \label{C:#1}  \ifShowLabels \TeXref{C:#1} \fi }%
        { \end{Cor} }
\newcommand{\cor}[1]{\begin{corol}{#1} \sl }
\newcommand{\ecor}{\end{corol}}

\newenvironment{defeni}[1]%
        { \begin{Def} \label{D:#1}  \ifShowLabels \TeXref{D:#1} \fi }%
        { \end{Def} }
\newcommand{\defe}[1]{\begin{defeni}{#1} \sl }
\newcommand{\edefe}{\end{defeni}}

\newenvironment{remark}[1]%
        { \begin{Rem} \label{R:#1}  \ifShowLabels \TeXref{R:#1} \fi }%
        { \end{Rem} }
\newcommand{\rem}[1]{\begin{remark}{#1}}
\newcommand{\erem}{\end{remark}}

\newenvironment{conjec}[1]%
        { \begin{Conj} \label{Co:#1}  \ifShowLabels \TeXref{Co:#1} \fi }%
        { \end{Conj} }
\renewcommand{\conj}[1]{\begin{conjec}{#1} \sl }
\newcommand{\econj}{\end{conjec}}

\newcommand{\eq}[1]%
        { \ifShowLabels \TeXref{E:#1} \fi
           \begin{equation} \label{E:#1} }
\newcommand{\eeq}{ \end{equation} }

\newcommand{\prf}{ \begin{proof} }
\newcommand{\epr}{ \end{proof} }

\newcommand\alp{\alpha}         

\newcommand\gam{\gamma}         \newcommand\Gam{\Gamma}
\newcommand\del{\delta}         \newcommand\Del{\Delta}
\newcommand\eps{\varepsilon}

\newcommand\kap{\kappa}
\newcommand\lam{\lambda}                \newcommand\Lam{\Lambda}
\newcommand\sig{\sigma}         

\newcommand\ome{\omega}         \newcommand\Ome{\Omega}

\newcommand\calA{{\mathcal{A}}}
\newcommand\calB{{\mathcal{B}}}
\newcommand\calC{{\mathcal{C}}}
\newcommand\calD{{\mathcal{D}}}

\newcommand\calF{{\mathcal{F}}}

\newcommand\calH{{\mathcal{H}}}

\newcommand\calK{{\mathcal{K}}}
\newcommand\calL{{\mathcal{L}}}
\newcommand\calM{{\mathcal{M}}}

\newcommand\calO{{\mathcal{O}}}
\newcommand\calP{{\mathcal{P}}}

\newcommand\calS{{\mathcal{S}}}

\newcommand\calZ{{\mathcal{Z}}}

            \newcommand\bfB{{\mathbf B}}
            \newcommand\bfC{{\mathbf C}}

\newcommand\bff{{\mathbf f}}            
            \newcommand\bfG{{\mathbf G}}
            
\newcommand\bfi{{\mathbf i}}

            \newcommand\bfL{{\mathbf L}}
\newcommand\bfm{{\mathbf m}}            \newcommand\bfM{{\mathbf M}}
            \newcommand\bfN{{\mathbf N}}
            \newcommand\bfO{{\mathbf O}}
\newcommand\bfp{{\mathbf p}}            \newcommand\bfP{{\mathbf P}}
\newcommand\bfq{{\mathbf q}}            
\newcommand\bfr{{\mathbf r}}            
\newcommand\bfs{{\mathbf s}}            
\newcommand\bft{{\mathbf t}}            \newcommand\bfT{{\mathbf T}}
\newcommand\bfu{{\mathbf u}}            \newcommand\bfU{{\mathbf U}}
            \newcommand\bfV{{\mathbf V}}
            
\newcommand\bfx{{\mathbf x}}            \newcommand\bfX{{\mathbf X}}
            \newcommand\bfY{{\mathbf Y}}
            \newcommand\bfZ{{\mathbf Z}}

\newcommand\QQ{\mathbb{Q}}

\newcommand\RR{\mathbb{R}}

\renewcommand\AA{\mathbb{A}}

\newcommand\GG{\mathbb{G}}

\newcommand\ZZ{\mathbb{Z}}

\newcommand\CC{\mathbb{C}}

 \newcommand\grg{{\mathfrak{g}}}

 \newcommand\grm{{\mathfrak{m}}}

 \newcommand\grp{{\mathfrak{p}}}

\newcommand\sdp{\times \hskip -0.3em {\raise 0.3ex
\hbox{$\scriptscriptstyle |$}}} 

\newcommand\Aut{\operatorname{Aut}}

\newcommand\End{\operatorname{End\,}}

\newcommand\gl{{\bf gl}}

\newcommand\Hom{\operatorname {Hom}}

\newcommand\Id{\operatorname{Id}}

\newcommand\Ind{\operatorname{Ind}}

\newcommand\Ker{\operatorname{Ker}}

\newcommand\RE{\operatorname{Re}}
\newcommand\Res{\operatorname{Res}}

\newcommand\Spec{\operatorname{Spec}}

\newcommand\supp{\operatorname{supp}}

\newcommand\tr{\operatorname{tr}}
\newcommand\Tr{\operatorname{Tr}}

\newcommand\oF{{\overline{F}}}

\newcommand\oT{{\overline{T}}}

\newcommand\oV{{\overline{V}}}

\newcommand\oX{{\overline{X}}}
\newcommand\oy{{\overline{y}}}

\newcommand\oPhi{{\overline{\Phi}}}

\newcommand\oome{\overline{\ome}}

\newcommand\tilf{{\widetilde{f}}}

\newcommand\tilY{{\widetilde{Y}}}

\newcommand\tilome{{\widetilde{\ome}}}

\newcommand\tilmu{{\widetilde{\mu}}}

\renewcommand\tilome{{\widetilde{\ome}}}

\newcommand\tilPhi{{\widetilde{\Phi}}}

\newcommand\tiltet{{\widetilde{\theta}}}

\newcommand\x{\times}
\newcommand\ten{\otimes}

\newcommand{\ra}{\rangle}
\newcommand{\la}{\langle}

\newcommand\qlb{{\overline \QQ}_l}
\newcommand\fr{\text{Fr}}
\newcommand\frw{\fr_w}
\newcommand\gln{{\bf GL}(n)}
\newcommand\wt{\widetilde}
\newcommand\chl{\calK_{\calL}}
\newcommand\ch{\text{ch}}

\newcommand{\bPhi}{\Phi\hskip-7.5pt \Phi}
\renewcommand{\oPhi}{{\overline \Phi}}
\renewcommand{\gl}{{\bf GL}}

\renewcommand{\Id}{\text{Id}}
\newcommand{\Irr}{\operatorname{Irr}}
\newcommand\gd{{\bfG}^{\vee}}
\newcommand\cic{C^{\infty}_c}
\renewcommand\Spec{\operatorname{Spec}}

\newcommand\tbY{{\widetilde \bfY}}
\newcommand\tbf{{\widetilde \bff}}
\newcommand\tbp{{\widetilde \bfp}}
\newcommand\obY{{\overline \bfY}}
\newcommand\obX{{\overline \bfX}}
\newcommand\obf{{\overline \bff}}
\newcommand\obp{{\overline \bfp}}
\renewcommand\oome{{\overline \ome}}
\newcommand\tbphi{{\widetilde \bPhi}}
\newcommand\obPhi{{\overline \bPhi}}
\newcommand\tbG{{\widetilde \bfG}}
\newcommand\obT{{\overline \bfT}_\rho}

\renewcommand\top{\text{top}}
\newcommand\obg{{\overline \bfG}_\rho}

\newcommand\td{\bfT^{\vee}}
\newcommand\gaf{\bfG(\AA_K)}
\renewcommand\cic{C^{\infty}_c}
\newcommand\phpsi{\bPhi_\psi}
\newcommand\phro{\bPhi_{\rho,\psi}}

\renewcommand\o{\omega}
\begin{document}
\title{$\gam$-functions of representations and lifting}
\author{A.~Braverman and D.~Kazhdan
(with an Appendix by V.~Vologodsky)}
\begin{abstract}
Let $F$ be a local non-archimedian field and let $\psi:F\to \CC^*$ be a
non-trivial additive character of $F$. To this data one associates a
meromorphic
function $\gam_\psi:\pi\mapsto\gam_\psi(\pi)$ on the set of irreducible
representations of the group $\gl(n,F)$ in the following way. Consider the
invariant distribution $\Phi_\psi:=\psi(\tr(g))|\det(g)|^n |dg|$ on $GL(n,F)$, where
$|dg|$ denotes a Haar distribution on $\gl(n,F)$. Although the
support of $\Phi_\psi$ is not compact, it is well-known that for generic
irreducible representation $\pi$ of $\gl(n,F)$ the action of
$\Phi_\psi$ in $\pi$ is well-defined
and thus it defines a number $\gam_\psi(\pi)$. These gamma-functions
and the associated L-functions were studied by J.~Tate for $n=1$ and by
R.~Godement and H.~Jacquet for arbitrary $n$.

Let now $G$ be the group of points of an arbitrary quasi-split reductive
algebraic group over $F$ and let $\bfG^{\vee}$ be the Langlands dual group. Let
also $\rho:\bfG^{\vee}\to \gl(n,\CC)$ be a finite-dimensional representation of
$\bfG^{\vee}$. The local Langlands conjectures predict
the existence of a natural
map $l_\rho:\Irr(G)\to\Irr(\gl(n,F))$. Assuming that a
certain
technical condition on $\rho$ (guaranteeing that the image of $l_\rho$ does not
lie in the singular set of $\gam_\psi$) is satisfied, one can consider the
meromorphic function $\gam_{\psi,\rho}$ on $\Irr(G)$, setting
$\gam_{\psi,\rho}(\pi)=\gam_\psi(l_\rho(\pi))$.

The main purpose of this paper is to propose a general framework for an
explicit construction of the functions $\gam_{\psi,\rho}$.
Namely, we propose a
conjectural scheme for constructing an invariant distribution
$\Phi_{\psi,\rho}$
on $G$, whose action in every $\pi\in \Irr(G)$ is given by multiplication by
$\gam_{\psi,\rho}$. Surprisingly, this turns out to be connected
with
certain geometric analogs
of M.Kashiwara's crystals (cf. \cite{BerK}).

We work out in detail several examples. As a byproduct we obtain a conjectural
formula for the lifting $l_\rho(\theta)$ where $\theta$ is a character
of an arbitrary maximal torus in $\gl(n,F)$. In some of these examples we also
give a definition of the corresponding local L-function and state a conjectural
$\rho$-analogue
of the Poisson summation formula for Fourier transform. This conjecture
implies that the  corresponding global $L$-function has meromorphic
continuation and satisfies a functional equation.
\end{abstract}

\maketitle

\sec{int}{Introduction and statement of the results}
\ssec{}{The Langlands program and the lifting}
Let $F$ be either local or finite field, $\oF$ -- its
separable closure. We will denote algebraic varieties
over $F$ and maps between them by boldface letters.
The corresponding ordinary letters will denote the
associated sets of $F$-points (and the induced maps
between them).

Let  $\bfG$ be a connected
reductive
algebraic group over $F$, $G=\bfG(F)$ its group of $F$-points,
$\Irr(G)$ the set of isomorphism classes of complex irreducible
representations of
$G$. For simplicity we shall assume in the introduction that $\bfG$ is
split, although
most of our results generalize easily to quasi-split groups. To
$\bfG$ one associates a connected complex algebraic group $\gd$,
called the Langlands dual
group.

\medskip
\noindent
{\bf Example.} If $\bfG=\gl(n)$ then $\gd=\gl(n,\CC)$.

\medskip
\noindent
The local Langlands conjectures can now be summarized as follows.
To the field
$F$ one associates a pro-algebraic group $\grg_F$ over $\CC$, which is very closely
related
with the Galois group $\text{Gal}(\oF/F)$ (in the case when $F$ is a
local
field
it is called the Weil-Deligne group).

The Langlands conjecture says that

1) There exists a finite-to-one map from $\Irr(G)$ to the set of all
$\gd$-orbits on  $\Hom(\grg_F,\gd)$.

2) If $ G=\gl(n,F)$ then the above map is a bijection.

\noindent
{\it Remark.} This conjecture is known for finite field (due to
Lusztig). It is  known for local fields when $\bf G$ is a torus (local class
field
theory) and for for $\gl(n)$ (cf. \cite{LRS},
\cite{Har} and \cite{HT}).

\medskip

Let now $\rho:\gd\to\gl(n,\CC)$ be an algebraic  representation of
$\gd$. Then
the above conjecture implies that there exists a map
$l_\rho:\Irr(G)\to\Irr(\gl(n,F))$.

Our main task is to try to propose a conjecture for an  explicit
description of the map $l_\rho$ and to formulate certain
$\rho$-analogue of the Poisson summation formula.

In the rest of this paper, until \refs{finite} we assume that $F$ is a
local field. Moreover, we assume that the characteristic of $F$ is
good
with respect to $\bfG$ and  that the Lie algebra $\grg$ of $\bfG$
possesses a non-degenerate invariant form. We denote by $\calH(G)$ the
Hecke algebra of $G$. By definition when $F$ is non-archimedian
(resp. archimedian)
this is the algebra of locally
constant (resp. $C^{\infty}$)
compactly supported distributions on $G$. We shall choose a
Haar
measure $|dg|$ and thus identify $\calH(G)$ with the space $\cic(G)$ of locally
constant (resp. $C^{\infty}$ when $F$ is archimedian)
compactly supported functions on $G$.
\ssec{int-gam}{$\gam$-functions}
We will not try to construct explicitly a representation $l_\rho
(\pi),\pi \in \Irr(G)$ but give an indirect description of  the lifting
$l_\rho$.

For this purpose we define a complex-valued function  $\gam$ on
$\Irr(\gl(n,F))$.

Choose (once and for all) a non-trivial additive character $\psi:F\to
\CC^*$.
Let $|dg|$ denote the unique Haar measure on $G=\gl(n,F)$ such that
the Fourier transform $\phi\mapsto\calF(\phi)$ defined by the
formula
\eq{}
\calF(\phi)(y)=\int\limits_G \phi(x)\psi(\tr(xy))|\det(x)|^n |dx|
\end{equation}
is unitary.

Consider the distribution $\Phi_\psi(g):=\psi(\tr(g))|\det(g)|^n |dg|$ on
$\gl(n,F)$.
Clearly, the distribution $\Phi=\Phi_\psi$
is invariant under the adjoint action. Let now $(\pi,V)$ be an irreducible
representation of $G$. Consider the integral
\eq{gam-act}
\pi(\Phi)=\int\limits_{G} \pi(g)\Phi(g)
\end{equation}
One can show (cf. \cite{GJ}) that the integral
\refe{gam-act} is convergent (although not absolutely convergent)
for
generic $\pi\in \Irr(G)$. Therefore, for generic $\pi$ the operator
$\pi(\Phi)$ is well-defined and it is equal to multiplication
by a scalar $\gam(\pi)$. One can consider  $\gam$ as a meromorphic function
on $\Irr(G)$ (cf. \refss{rat} for the definition of this notion).

Let now $G$ and $\rho$ be as above and
assume that we know the lifting map $l_\rho$. Then we define a
function $\gam_{\psi,\rho}=\gam_\rho$ on the set $\Irr(G)$ by
$\gam_{\rho}=l_\rho ^{*}(\gam)$. In other words any
irreducible representation $\pi$ of $G$ we define

\eq{gamro}
\gam_{\rho}(\pi)=\gam(l_\rho(\pi))
\end{equation}

The function $\gam_\rho$ is a well-defined meromorphic function
on $\Irr(G)$ provided that $l_\rho(\pi)$ does not lie in the singular
set
of $\gam$ for generic $\pi\in\Irr(G)$. In order to guarantee this we
shall always (except for \refs{finite}) assume that $\rho$ satisfies
the following condition:
there exists a cocharacter $\sig:\GG_m\to Z(\gd)$ of the center $Z(\gd)$
of the group $\gd$
such that $\rho \circ\sig= \Id$. Note that $\sig$ can be
regarded as a character of $\bfG$.

It is easy to show (assuming that $l_\rho$ satisfies certain natural
properties) that there exists an ad-invariant distribution
$\Phi_\rho=\Phi_{\rho,G}$
on $G$ such
that
\eq{gam}
\gam_\rho(\pi)=\pi(\Phi_\rho)\cdot \Id
\end{equation}

We see from \refe{gam} that a lifting $l_\rho$ of  a representation
$\rho$ of the dual group $\gd$ determines an
invariant distribution $\Phi_\rho$ on $G$.
Since the map $l_\rho$ collapses the $L$-packets we assume that
the function $\Phi_{\rho,G}$  is stable.
Our main purpose is to give a
conjectural description of the function  $\Phi_{\rho,G}$. In general we
can only propose a framework for such a  description. Sometimes  we can
make our
suggestion precise and  in some simplest cases  when the lifting is
known  we can check that our
definition is correct.
\ssec{}{The Fourier transform $\calF_\rho$}

Assume now that we know the lifting $l_\rho$ and therefore can construct
 the distribution $\Phi_{\rho}$ on
$G=\bfG(F)$ .
In \refs{poisson} we define certain Fourier-type transform operator
$\calF_\rho$ acting in the space of functions on $G$. Namely, for
$\phi\in\calH(G)$ we define
\eq{}
\calF_\rho(\phi)=|\sig|^{-l-1}(\Phi_\rho*{^\iota\phi})
\end{equation}
where $^\iota\phi(g)=\phi(g^{-1})$ and $l$ is the semi-simple rank of
$\bfG$.

 In the case when $\bfG=\gln$ and $\rho$ is the
standard
representation of $\gd\simeq\gln$ the operator $\calF_\rho$ coincides
with the usual Fourier transform in the space of functions on the space
$\bfM_n(F)$ of $n\x n$-matrices over $F$.

We conjecture that $\calF_\rho$ extends to a unitary operator acting
in the space $L^2_\rho(G)=L^2(G,|\sig|^{l+1}|dg|)$ where $|dg|$ is a Haar
measure on $G$ and that in the case when
$\rho$ is non-singular the space $\calS_\rho(G)$ is stable under
$L^2_\rho(G)$.

Note that in \refs{action} we give a definition of $\Phi_\rho$ ( and thus of
$\calF_\rho$) in certain cases where
$l_\rho$ is not known.

\ssec{lfunctions}{Local $L$-functions}

Until now we have discussed only the local lifting. But the only known
formulation of the {\it local lifting conjecture} is to
interpret it as a local
part of the {\it global lifting conjecture}. There are two approaches to a
proof of   the {\it global lifting conjecture}. The first approach is based
on the
Trace formula and the second one on the study of $L$-functions. So it is
natural to try to give a direct definition of the $L$-function
$L(l_\rho (\pi),s)$ for $\pi \in \Irr(G)$ in terms of $\rho$ and $\pi$.
Assume that we know the distribution $\Phi_\rho$. In such a case 
we present a conjecture  for a  direct definition of the $L$-function
$L(l_\rho (\pi),s)$. To any 
representation   $\rho$ of $\gd$ we associate
a certain subspace $\calS_\rho$ of the
space of functions on $G$ such that  $\calS_\rho$ contains the space
$\calS(G)$ of smooth compactly supported functions on $G$. For any
representation $\pi \in \Irr(G)$ we
define the function $L(l_\rho (\pi),s)$ as the common denominator of
rational operator-valued functions $\pi _s (\phi)$ where $\phi\in \calS_\rho$
and  $\pi _s (f):=\int _{g\in G}\phi(g)|\sig (g)|^s \pi (g)|dg|$. In the case
when $\bfG=\gln$ and $\rho$ is the standard
representation of  $\gd\simeq\gln$ the space  $\calS_\rho$ coincides
with the space of Schwartz-Bruhat functions on the space $M_n$ of $n\times
n$-matrices and our definition of the  $L$-function coincides with the
definition from  \cite{GJ}.

\ssec{}{Conjectural applications to automorphic $L$-functions}Let now $K$ denote
a global field and let also $\AA_K$ be the corresponding adele ring.
For a place $\grp$ of $K$ we denote by $K_\grp$ the corresponding
local
completion of $K$.
In \cite{GJ} R.~Godement and H.~Jacquet
defined the L-function $L(\pi,s)$ for every automorphic
representation
$\pi$ of $\gl(n,\AA_K)$. This $L$-function is defined as the product of
the corresponding local $L$-functions over all places of $K$.
It is shown in \cite{GJ} that $L(\pi,s)$ is meromorphic and has a functional
equation. The main tools in the construction and the proof are the
space $\calS(\bfM_n(F))$ of Schwartz-Bruhat functions on $n\x n$-matrices,
the Fourier transform acting in this space and the  Poisson summation
formula .

Assume now that we are given a group $\bfG$ and  a
representation $\rho$
of $\gd$. R.~Langlands (cf. \cite{Lan}) conjectured that one could define
an L-function $L_\rho(\pi,s)$ attached to every automorphic
representation $\pi$
of $\bfG(\AA_K)$ as a product of local factors. We conjecture that these
local factors are equal to ones defined  in \refss{lfunctions}. In
\refs{local} we define a
global version $\calS_\rho(\bfG(\AA_K))$ and define a global
$\rho$-analogue  $\calF_\rho$ of the Fourier transform. By the
definition the space  $\calS_\rho(\bfG(\AA_K))$ is the span of functions
$\rho$ which are products
$\phi=\bigotimes \phi_{\grp}$ of local factors. We
conjecture the
existence of an  $\calF_\rho$-invariant distribution $\eps _\rho $ on
$\calS_\rho(\bfG(\AA_K))$ such that

1)
\eq{}
\eps_\rho (\phi)=\sum\limits_{g\in \bfG(K)}\phi(g)
\end{equation}
 if some local factor $\phi _{\grp}$  has compact support on $G_{\grp}$

2)For any $\phi\in\calS_\rho(\gaf)$
\eq{}
\eps_\rho(\phi)=\eps_\rho(\calF_\rho(\phi))
\end{equation}

Existence of such $\eps$ can be thought of as  a $\rho$-analogue
of the Poisson summation formula. In the case when  $\bfG=\gln$ and
$\rho$ is the standard
representation of $\gd\simeq\gln$ our conjecture specializes to the
usual  Poisson summation formula. The validity of our  $\rho$-analogue
of the Poisson summation formula would imply that the automorphic
$L$-function $L_\rho(\pi,s)$ has meromorphic continuation with only a
finite number of poles and satisfies a functional equation.
\ssec{litor}{Application to lifting}In \refs{action} we give a  construction
of the distribution $\Phi_\rho$ in the case when $\bfG=\gln\x \gl(m)$
for arbitrary $n$ and $m$ and when $\rho$ is equal to the tensor product
of the
standard representations of $\gl(n,\CC)$ and $\gl(m,\CC)$. As a
byproduct we get
a conjectural formula for $l_\rho(\theta)$ where $\theta$ is a character
of an
arbitrary maximal torus in $\gl(n,F)$.
\ssec{}{Idea of the construction of $\Phi_\rho$}
In \refs{torus} we check that in the case when
$\bfG$ is a split torus  $\bfT$ the lifting and the distribution
$\Phi_{\rho,T}$ on $T$ are well
defined for any representation $\rho_T$ of ${\bf T}^\vee$, satisfying the
condition from \refss{int-gam}.  Moreover,
one can construct an
algebraic variety $\bfY_{\rho,\bfT}$, a map ${\bf p_T:Y_{\rho,T}
\rightarrow T}$, a function $\bff_\bfT$ on
$\bfY_{\rho,\bfT}$ and a top-form $\ome_\bfT \in  \Gam(\Omega ^{top},
\bfY_{\rho,\bfT} )$
such that the distribution $\Phi
_{\rho,T}$ is equal to the push-forward $(p_T)_! (\psi (f_T)|\o_T |)$
where $|\o_T|$
is the measure on $ Y_{\rho,T}$ associated with the form $\ome_T$
(cf. \cite{weil}).
We
say that $\bPhi_{\rho,\bfT} =( \bfY_{\rho,\bfT} ,\bff_{\rho,\bfT} ,\ome_{\rho,\bfT} )$
is an {\it algebraic-geometric distribution representing} the distribution
$\Phi _{\rho ,T}$ and say that  the distribution$\Phi _{\rho ,T}$ is a
{\it materialization} of an  algebraic-geometric distribution
$\bPhi_{\rho,\bfT}$.

Let now $\bfG$ be an arbitrary reductive group. We conjecture that the
$Ad$-invariant distribution  $\Phi _{\rho}$ is $stable$ and therefore comes
from a distribution   $\tilde {\Phi} _{\rho}$ on the set ${\bf G/ Ad} (F)$ of
$F$-rational
points of the geometric quotient  $\bf G/ Ad =\Spec(\calO(\bfG))^{\bfG}$ of
$\bfG$ by the adjoint action. Since
 ${\bf G/ Ad=T}/W$  where  $\bf T \subset \bf G$ is a maximal split
torus and  $W$ is the Weyl group  of  $\bf G$ we can consider  $\tilde
{\Phi} _{\rho}$ as a distribution on the set $({\bfT}/W)(F)$. We conjecture
that the distribution  $\tilde {\Phi} _{\rho}$  is a
{\it materialization} of  an algebraic-geometric distribution
$\bPhi_\rho$ which
is a ``descent'' of  $\bPhi_{\rho_\bfT,\bfT}$ where $\rho_\bfT$ is the restriction
of $\rho$ on  ${\bf
T}^\vee$. More precisely we conjecture that there exists an action of the
Weyl group $W$ on  $\bPhi_{\rho,T}$  such that the distribution $\tilde
{\Phi} _{\rho}$ on
the set  ${\bf T}/W(F)$  is the
{\it materialization} of the   algebraic-geometric
distribution  $\bPhi_\rho$ on ${\bf T}/W$ which is a "descent" of  $\bPhi_{\rho_\bfT,\bfT}$
(we shall explain the ``descent'' procedure carefully
in \refs{algint}).

In \refs{action} we give explicit formulas for this action in a number of cases.
\ssec{}{}This paper is organized as follows. In \refs{lifting} we
formulate the Langlands lifting conjecture and discuss several
examples. In \refs{gamma} we define $\gam$-functions and formulate
their conjectural properties.
In \refs{torus} we give an explicit construction of
the lifting
for the case when our
reductive group is a split torus
$\bf T$ and we write an explicit formula for the corresponding
distribution $\Phi_{\rho,T}$.
\refs{local} and \refs{poisson} are
devoted respectively
to the definition of the Schwartz space $\calS_\rho$ and an analogue
of the Poisson summation formula for the operator $\calF_\rho$.
In \refs{algint} we develop the notion of an algebraic-geometric
distribution. Using this notion we reduce (conjecturally) the problem of constructing
the distribution
$\Phi_\rho$ to a purely algebraic question.

\refs{action} is devoted to the discussion of this question
in the following case. Let $m,n$ be two positive integers.
Define the group
\eq{}
\bfG(m,n)=\{ (A,B)\in\gl(m)\x \gln|\ \det(A)=\det(B)\}.
\end{equation}
The dual group $\bfG(m,n)^{\vee}$ is the quotient
of $\gl(m,\CC)\x \gl(n,\CC)$ by the subgroup consisting of all pairs
of matrices $(t\Id_m,t^{-1}\Id_n)$ for $t\in\CC^*$. Hence the tensor
product $\rho_m\ten\rho_n$ of the standard representations of
$\gl(m,\CC)$ and $\gl(n,\CC)$ descends to a representation
$\rho$ of $\bfG(m,n)^{\vee}$. The corresponding $\gam$-function
is essentially constructed in \cite{JPSS}.

In \refs{action}
we give a construction of $\Phi_\rho$ in this case.
Using \cite{BerK}
one can check that our definition gives rise to (almost) the same
$\gam$-function as the one defined in \cite{JPSS}. As an application
we give conjectural formulas for the lifting problem discussed in
\refss{litor}.

In \refs{finite} we explain how to construct $\Phi_\rho$ in
the case when the field $F$ is finite. More precisely, in this case we
construct a perverse sheaf $\bPhi_\rho$ on $\bfG$ such that
$\Phi_\rho$
is obtained from it by taking traces of the Frobenius morphism in the
fibers.

Finally, the Appendix (by V.~Vologodsky) contains a proof of  some
result about non-archimedian oscillating integrals used throughout the
paper.
\ssec{notations}{Notations}Everywhere, except for section \refs{finite},
$F$ denotes a local
field. If $F$ is non-archimedian then we denote by $\mathcal O _F\subset
F$ its ring of integers and  by $q$ the number of elements in the
residue field of $F$.

If  $X$ is either a totally disconnected topological space or a
$C^{\infty}$-manifold we denote by $C(X)$ the space of $\CC$-valued continuous functions on $X$.
We denote by $C^{\infty}_c(X)$ the space of all compactly supported $\CC$-valued functions
on $X$ which are:

1) locally constant if $X$ is totally disconnected;

2) $C^{\infty}$ if $X$ is a $C^{\infty}$-manifold.

We shall denote algebraic varieties over $F$ by boldface letters (e.g. $\bfG$,$\bfX$...).
The corresponding ordinary letters ($G$, $X$,...) will denote the corresponding
sets of $F$-points.

Let $\bfX$ be a smooth algebraic variety over $F$ and
$\ome\in\Gam(\bfX,\Ome^{top}(\bfX))$, where $\Ome^{top}$ denotes the sheaf
of differential forms of degree $\dim(\bfX)$ on $\bfX$. According to
\cite{weil} to $\ome$ one can associate a distribution $|\ome|$.

If $\bfG$ is a topological group (resp. a Lie group) then we denote
by $\Irr(G)$ the set of isomorphism classes of irreducible algebraic
(cf. \cite{BZ1}) representations of $G$ (resp. the set of equivalence
classes of continuous irreducible
Banach representations of $G$ with respect to infinitesimal
equivalence -- cf. \cite{Wa}).

\ssec{}{Acknowledgments}
We are grateful A.~Berenstein,
J.~Bernstein, V.~Drinfeld, D.~Gaitsgory, M.~Harris,
M.~Kontsevich, S.~Rallis and V.~Vologodsky for very helpful and interesting
discussions
on the subject. We are also grateful to Y.~Flicker for numerous
remarks about this paper.

\sec{lifting}{Lifting}

\ssec{unramified}{Unramified lifting}
Let $F$ be a non-archimedian local field and $\bfG$ a quasi-split
group over $\calO_F$, which can be
split over an unramified extension of $F$ of degree $k$ (we
assume that $k$ is minimal with this property). Then the cyclic
group $\ZZ/k\ZZ$ acts naturally on $\gd$. We denote the action
of the generator of $\ZZ/k\ZZ$ by $g\mapsto \kap(g)$.

The group $\gd$ acts on itself by $g:x\mapsto gx\kap(g)^{-1}$.
Let $S_\bfG$ denote the set of closed orbits on $\gd$ with respect
to the above action.

Let $\Irr_{un}(G)$ denote the set of isomorphism classes of irreducible
unramified
representations of $G$, i.e. the set of irreducible representations
of $G$ which have a non-zero $\bfG(\calO_F)$-invariant
vector. Then one has the natural identification
$\Irr_{un}(G)\simeq S_\bfG$.

Let $\rho:\gd\rtimes\ZZ/k\ZZ\to \gln$ be a representation. The map
$x\mapsto (x,\kap)$ gives rise to a well-defined map from
$S_\bfG$ to the set of semi-simple conjugacy classes in $\gd\rtimes\ZZ/k\ZZ$.
If we compose this map with $\rho$ we get a map from
$S_\bfG$ to $S_{\gln}$.

Consider now
the composite map
\eq{}
\Irr_{un}(G)\to S_\bfG\to S_{\gln}\to\Irr_{un}(\gl(n,F)).
\end{equation}
We will denote this map by $l_{\rho,un}$ and call it
{\it the unramified lifting}.
\ssec{glob-lift}{The global lifting}
Let $K$ be a global field, $\AA_K$ -- its ring of adeles. We denote
by $\calP(K)$ the set of places of $K$. For every $\grp\in\calP(K)$ we
let $K_\grp$ be the corresponding local completion of $K$.

Let now $\bfG$ be a quasi-split reductive group over $K$, which
splits over a finite separable extension of $K$ with Galois group
$\Gam$.
As is well-known any irreducible representation $\pi$ of the
group $\bfG(\AA_K)$ can be uniquely written
as a restricted tensor product
$\ten_{\grp\in\calP(K)}\pi_\grp$ where
$\pi_\grp\in\Irr(\bfG(K_\grp))$ and
$\pi_\grp\in\Irr_{un}(\bfG(K_\grp))$
for almost all $\grp$.

We denote by
$\Aut(\bfG,K)$ the set of irreducible automorphic representations of
$\bfG(\AA_K)$.

The group $\Gam$ acts naturally on $\gd$. Hence we can consider the
semidirect
product $\gd\rtimes \Gam$. For almost all $\grp\in \calP(K)$ the
group $\bfG$ can be split over an unramified extension of $K_\grp$
of degree $k_\grp$ and in this case the group $\gd\rtimes
\ZZ/k_\grp\ZZ$ is naturally embedded in  $\gd\rtimes \Gam$.

Hence every representation $\rho:\gd\rtimes\Gam\to\gln$
defines a representation of $\gd\rtimes \ZZ/k_\grp\ZZ$ for almost
all $\grp$.
\defe{}Let $\pi=\ten \pi_\grp$ be an automorphic representation of
$\bfG(\AA_K)$. The {\it global lifting} $l_\rho(\pi)$ is an
automorphic representation $l_\rho(\pi)=\ten  l_\rho(\pi)_\grp$
of the group $\gl(n,\AA_K)$ such that for every place $\grp\in\calP(K)$
such that $\pi_\grp$ is unramified one has $l_\rho(\pi)_\grp=
l_{\rho,un}(\pi_\grp)$ where $l_{\rho,un}$ is as in \refss{unramified} .
\edefe

\noindent
{\it Remark.} It follows from the strong multiplicity one theorem
for $\gl(n)$ (cf. \cite{PS}) that if $l_\rho(\pi)$ exists then it is unique.

\conj{lifting}
\begin{enumerate}
\item
For every $\bfG$, $\pi$  and $\rho$ as above the
lifting $l_\rho(\pi)$ exists.
\item
Let $\bfG$ be a connected quasi-split reductive group over a local field
$F$, which
splits over a normal separable extension $L/F$ such that
$\text{Gal}(L/F)=\Gam$. Then for
every representation $\rho:\gd\rtimes \Gam\to\gl(n,\CC)$
there exists a map $l_\rho ^F:\Irr(G)\to\Irr(\gl(n,F))$ such that for any
global field  field $K$, a place  $\grp\in\calP(\calK)$ such that
$F=K_\grp$, a quasi-split group
$\bfG '$ over $K$ as in \refss{glob-lift} such that ${\bfG}'_\grp\simeq\bfG$
and  every $\pi\in\Aut(\bfG ',K),\pi=\ten \pi_\grp$, the $\grp$-th
component of the automorphic   representation $l_\rho(\pi)$ is equal to
$l_\rho ^F ( \pi_\grp )$.
\end{enumerate}
\econj

We say that $\rho$ is liftable if \refco{lifting} holds for $\rho$.

\noindent
{\it Remark.} It is easy to see that the local lifting $l_\rho$
is uniquely determined by the conditions of \refco{lifting}.


\ssec{lift-torus}{Example}One of the examples that we would like
to consider in this paper
is the following. Let $E/F$ be a separable extension of $F$ of degree $n$.
Let $\bfT_E=\Res_{E/F}\GG_{m,E}$ where $\Res_{E/F}$ denotes the
functor of restriction of scalars. Let $\Gam_F$ denote the absolute
Galois group of $F$.
$\Gam_F$ acts naturally on the set of all embeddings
of $E$ into $\oF$, which is a finite set with $n$ elements.
Hence we get a homomorphism $\alp_E:\Gam_F\to S_n$ which is defined
uniquely
up to $S_n$-conjugacy. Let $\Gam=\Im (\alp_E)$. By the definition $\Gam$
is a subgroup of $S_n$ which is also a quotient of $\Gam_F$.
Moreover, it is clear that the torus $\bfT_E$ splits over the
Galois extension
$L/F$ where $\text{Gal}(L/F)=\Gam$.

The dual torus $\td$ is isomorphic to $(\CC^*)^n$ and the group
$\Gam$ acts on it by the restriction of the standard action of $S_n$
on $(\CC^*)^n$ to $\Gam$. Hence the standard embeddings of
$(\CC^*)^n$ and $S_n$ to $\gl(n,\CC)$ give rise to a homomorphism
$\rho _{\bfT}:\td\rtimes \Gam\to\gl(n,\CC)$.
One of the purposes of this paper
is to give a conjectural construction of $l_\rho$ in this case.

\sec{gamma}{Gamma-functions and Bernstein's center}
\ssec{}{Bernstein's center}Let $\bfG$ be a reductive algebraic group
over
$F$, $G=\bfG(F)$. We denote by $\calM(G)$ the category of smooth
representations
of $G$.

Recall that the {\it Bernstein center} $\calZ(G)$ of $G$ is the algebra
 of endomorphisms of the identity functor on the
category $\calM(G)$. It was shown by Bernstein that $\calZ(G)$ is isomorphic to
the direct product of algebras of the form $\calO(\Omega)$ where $\Omega$ is a
finite-dimensional irreducible complex algebraic variety. We shall denote the set
of all $\Ome$ as above by $\calC$. By $\Spec(\calZ(G))$ we
shall
mean the disjoint union of all $\Ome\in\calC$. It is shown in
\cite{BD}
that the natural map
$\Irr(G)\to \Spec(\calZ(G))$ sending every irreducible representation to its
infinitesimal character is finite-to-one and moreover is  generically
one-to-one.

The following result (also due to J.~Bernstein) gives a different
interpretation
of $\calZ(G)$. Let $\Phi$ be a distribution on
$G$. We say that $\Phi$ is essentially compact if for every function
$\phi\in \calH(G)$ one has $\Phi*\phi\in \calH(G)$. Given two essentially
compact distributions $\Phi_1,\Phi_2$ one defines their convolution
$\Phi_1*\Phi_2$ by setting
\eq{}
\Phi_1*\Phi_2(\phi)=\Phi_1(\Phi_2*\phi)
\end{equation}
It is easy to see that $\Phi_1*\Phi_2$ is again an essentially compact
distribution.
\lem{bern-center}
The algebra $\calZ(G)$ is naturally isomorphic to the algebra of invariant
essentially compact distributions.
\elem
\ssec{variant}{Example}Let $\bfG=\gln$, $G=\gl(n,F)$. Fix any $a\in F^*$ and a non-trivial
additive character $\psi:F\to\CC^*$. Let
$G_a=\{g\in G|\det(g)=a\}$. Then a choice of a Haar measure on
${\bf SL}(n,F)$ defines a measure on $G_a$ for
every $a\in F^*$ (since $G_a$ is a principal homogeneous
space over ${\bf SL}(n,F)$). We denote this measure by $d_a g$.
Define a distribution $\Phi_a$ on $G$ supported on $G_a$ by
\eq{}
\Phi_a(\phi)=\int\limits_{G_a}\phi(g) \psi (\tr (g)) d_a g
\end{equation}
for any $\phi\in C^{\infty}_c(G)$.

One can verify that $\Phi_a$ is essentially compact.
\ssec{rat}{Rational functions on $\Irr(G)$}Let $\calK(G)$ denote the
``field of fractions'' of $\calZ(G)$, i.e. the associative algebra
$\prod\calK(\Ome)$, where $\Ome$ runs over connected components of
$\Spec(\calZ(G))$ as above and $\calK(\Ome)$ is the field of fractions
of $\calO(\Ome)$. Since  generically
$\Spec(\calZ(G))$ and $\Irr(G)$ are identified we shall refer to
elements of $\calK(G)$ as {\it rational
functions} on $\Irr(G)$.
\ssec{}{$\sig$-regular central elements}For an associative algebra
$\calA$ we let $\calM(G,\calA)$ denote the category of $G$-modules
endowed with an action of $\calA$ by endomorphisms of the
$G$-module structure (we call it the category of
$(G,\calA)$-modules).

Let $\calA=\CC[t,t^{-1}], \calB=\CC((t))$.
Assume that we are given an $F$-rational character
$\sig:\bfG\to\GG_m$. Then we can define an action
of $G$ on $\calA$ and $\calB$ by requiring that
\eq{}
g\cdot t=|\sig(g)| t
\end{equation}
for every $g\in G$.

Define the functor $F_\calA:\calM(G,\calA)\to \calM(G,\calB)$ such
that
$F_\calA:V\mapsto V\underset{\calA}\otimes \calB$.
Set $\calZ_\sigma(G)=\End F_\calA$ (the $\End$ is taken in
the category of functors from $\calM(G,\calA)$ to $\calM(G,\calB)$).

We would like now to identify $\calZ_\sig(G)$ with a subspace of
$\calK(G)$. The assignment $\pi\mapsto\pi\ten |\sig|^{\text{log}_q(z)}$
gives rise to an action of $\CC^*$ on $\Irr(G)$ (here $z\in
\CC^*$). It is easy to see that this action descends to a well-defined
algebraic
action of $\CC^*$ on $\Spec(\calZ(G))$. For every $f\in \calK(G)$ we
can
write
\eq{}
f(z\cdot x)=\sum f_n(x)z^n
\end{equation}
for some $f_n\in\calK(G)$.
\lem{sreg}$\calZ_\sig(G)$ is naturally isomorphic to the space of all
$f\in \calK(G)$ such that $f_n\in\calZ(G)$ for every $n\in \ZZ$.
\elem
\refl{sreg} explains the term ``$\sig$-regular''.
\ssec{}{$\sig$-compact distributions}Let now $\Phi$ be a distribution
on $G$ and let $G_n=\sig^{-1}({\mathfrak \pi}^n\calO^*)$. $G_n$ is an open
subset of $G$. For every $n\in \ZZ$ we define a new distribution
$\Phi_n$ on $G$ by $\Phi_n=\chi_n\Phi$, where $\chi_n$ is the
characteristic
function of $G_n$.
\defe{scomp}
We say that $\Phi$ is {\it $\sig$-compact} if the following two
conditions are satisfied:
\begin{enumerate}
\item
$\Phi_n$ is essentially
compact for every $n\in \ZZ$.
\item For every $\phi\in\calH(G)$ there exists a rational
function $F_\phi:\CC\to \calH(G)$ such that
for every $g\in G$ one has
\eq{}
\sum\limits_{n\in \ZZ}z^n(\Phi_n*\phi)(g)=F_\phi(z)(g)
\end{equation}
\end{enumerate}
\edefe
The following lemma is straightforward from the definitions and
\refl{bern-center}.
\lem{}The space of all $\sig$-compact distributions is naturally
isomorphic
to $\calZ_\sig(G)$.
\elem

\ssec{ex-sigcomp}{An example}Let $G=\gl(n,F)$. Choose a non-trivial additive
character
$\psi:F\to \CC^*$ and set
$\Phi_\psi(g)=\psi(\tr(g))|\det(g)|^n |dg|$, where $|dg|$ is a Haar measure on
$G$.
Let also $\sig:\gln\to\GG_m$ be given by $\sig(g)=\det(g)$.
\prop{fgln}$\Phi_\psi$ is $\sig$-compact.
\eprop
\prf Let $d_a g=|\det|^n |dg|$ be the ``additive'' measure on $G$ and
let also $\calS(M_n)$ denote the space of Schwartz-Bruhat functions
on $n\x n$-matrices $M_n$ over $F$. We shall regard
$\calS(M_n)$ as a subspace of $C(G)$. Consider the Fourier
transform $\calF:\calS(M_n)\to \calS(M_n)$ given by
\eq{fourier}
\calF(\phi)(y)=\int\limits_{G}\phi(x)\psi(\tr(xy))d_a x
\end{equation}
It is easy to see that $\calF(\phi)=(\Phi_\psi *{^\iota
\phi})|\det|^{-n}$, where ${^\iota\phi}(x)=\phi(x^{-1})$.

Let us now show that conditions 1 and 2 above hold. Without loss
of generality we may assume that $\supp \phi\subset G_0$.
Then $q^{-n}\Phi_n*{^\iota\phi}=\calF(\phi)\chi_{-n}$. Since
$\calF(\phi)\in\calS(M_n)$ it follows that
$\calF(\phi)\chi_{-n}\in\calH(G)$.
The verification of the second conditions is left to the reader.
\epr

We shall denote by $\gam_\psi$ the rational function on
$\Irr(\gl(n,F))$ which corresponds to $\Phi_\psi$.
\ssec{fpro}{Lifting and the distributions $\Phi_{\psi,\rho}$}Let now
$\rho:\gd\to \gl(n,\CC)$  be a finite-dimensional representation
of the Langlands dual group $\gd$.
The local lifting conjecture predicts the existence of the natural
map $l_\rho:\Irr(G)\to \Irr(\gl(n,F))$. Assuming that
we know $l_\rho$ we can try to define a rational function
$\gam_{\psi,\rho}$
by setting $\gam_{\psi,\rho}(\pi)=\gam_\psi(l_\rho(\pi))$.
However, it might happen that the image of $l_\rho$ lies
entirely inside the singular set of $l_\rho$.
To guarantee that this does not happen we introduce the following
notion.
\defe{non-deg}
We say that $\rho$ is admissible if the following conditions hold.
\begin{enumerate}
\item
$\Ker(\rho)$ is connected.
\item
There exists a character $\sig:\bfG\to \GG_m$ defined over $F$ such
that
\eq{}
\rho\circ\sig=Id
\end{equation}
where we regard $\sig$ as a cocharacter of the center $\bfZ(\gd)$ of
$\gd$ and $I_n\in\gl(n,\CC)$ denotes the identity matrix.
\end{enumerate}
\edefe

It is easy to see that if we assume that $\rho$ is admissible and liftable
then $\gam_{\psi,\rho}$ satisfies the conditions of
\refl{sreg}.
Thus we can construct the corresponding $\sig$-compact
distribution $\Phi_{\psi,\rho,G}$. Our main task will be to try to
construct the distribution $\Phi_{\psi,\rho,G}$ explicitly  (we will sometimes drop
the indices $\psi$ and $G$, when it does not lead to a confusion).
\ssec{basex}{An example}One of the examples that we would like
to study in this paper is the following.
Fix a sequence $n_1,...,n_k$ of
natural numbers. Let $\bfG=\bfG(n_1,...,n_k)$ be the subgroup
of $\gl(n_1)\x ...\x \gl(n_k)$ consisting of $k$-tuples of non-degenerate
matrices $(A_1,...,A_k)$ such that $\det A_i=\det A_j$ for all
$1\leq i,j\leq k$. The dual group $\gd$ is the quotient of
$\gl(n_1)\x ...\x \gl(n_k)$ by the subgroup consisting of the
$k$-tuples $(t_1\Id_{n_1},t_2\Id_{n_2},...,t_k\Id_{n_k})$ with
$t_1...t_k=1$. Hence the tensor product representation
$\rho_{n_1}\ten ...\ten\rho_{n_k}$
of $\gl(n_1)\x ...\x \gl(n_k)$
factors through $\gd$ and we define $\rho=\rho_{n_1}\ten
...\rho_{n_k}$. We shall refer to $\rho$ as {\it standard
representation} of $\gd$. It is easy to see that  $\rho$ is an embedding
and that $\rho\circ\sig=Id$ where $\sig =det\circ \rho$.

\medskip
\noindent
{\it Example.} Suppose that $k=1$, $n_1=n$.
In this case $G=\gl(n,F)$
and $\Phi_{\psi,\rho}$ is the same as in \refss{ex-sigcomp}.
\ssec{induction}{Compatibility with induction}Let $\bfP\subset \bfG$ be a
parabolic subgroup of $\bfG$ defined over $F$, $\bfM$ -- the
corresponding Levi subgroup. Let $\del_\bfP$ be the determinant of the
action of $\bfM$ on the Lie algebra of the unipotent radical
$\bfU_\bfP$ of
$\bfP$.

For every $\phi\in\calH(G)$ we can construct a new function
$r_P(\phi)$ on $G/U_P$
\eq{}
r_P(\phi)(g)=\int_{U_P} \phi(gu)|du|
\end{equation}
Note that $M$ acts on $G/U_P$ on the right.

Let now $\rho$ be a representation of $\gd$ as above. Since
$\bfM^{\vee}$ is naturally a subgroup of $\gd$ we can regard
$\rho$ as a representation of $\bfM^{\vee}$. We conjecture that
the following property holds for the distributions $\Phi_{\psi,\rho,G}$ and
$\Phi_{\psi,\rho,M}$:

\medskip
\noindent
\conj{compatibility}
\eq{}
r_P(\Phi_{\psi,\rho,G}*\phi)=r_P(\phi)*(\Phi_{\psi,\rho,M}|\del_\bfP|^{-1/2}).
\end{equation}
\econj
\ssec{dir-sum}{}
Assume that $\rho=\rho_1\oplus\rho_2$ is a direct sum of two
admissible representations, such that the corresponding character
$\sig$ is the same for both representations.
The following result can be easily deduced from the definitions.
\lem{}
One has
\eq{dir-sum}
\Phi_{\psi,\rho}=\Phi_{\psi,\rho_1}*\Phi_{\psi,\rho_2} .
\end{equation}
\elem

\ssec{kaz-def}{$\gam$-functions determine the lifting}
Let ${\bfG=\gl(n-1)\times \gl(n)}$ and the character $\sigma$ be given by
$\sig (g_{n-1} ,g_n)=det (g_{n-1})$. In \cite{GeKa} Ch.7 a
rational function $\gam _{n-1,n}$ was defined on the subset of
non-degenerate representations of the group  $G$. One can show that
$\gam_{n-1,n}$ extends to a $\sig$-regular
central element for $G$. Hence there exists a $\sig$-regular distribution $\Phi_{n-1,n}$ on
the group  $\gl(n-1,F)\x\gl(n,F)$ such that
$\gam _{n-1,n}=\gam _{\Phi_{n-1,n}}$. In
\cite{JPSS} the analogous function   $\gam _{m,n}$ on $\Irr
(\gl(m,F)\x\gl(n,F))$  was
defined. It is also easy to see that it comes from  a $\sig$-regular
distribution $\Phi_{m,n}$ on the group  $\gl(m,F)\times\gl(n,F)$.

As was shown in \cite{GeKa} for any non-degenerate representation $\pi$ of
$\gl(n,F)$ the rational function  $\gam _{\pi}\doteq \gam
_{n-1,n}(\pi,\star )$ on $\Irr \gl(n-1,F)$
determines  uniquely the representation $\pi$. Moreover \cite{GeKa} contains
an explicit recipe for the construction of the representation $\pi$ in
terms of the function  $\gam _{\pi}$. It is clear from this recipe that we
have to know the  function  $\gam _{\pi}$ only up to multiplication by a
constant.
\ssec{}{}Let now $\bfG$ be arbitrary and
$\rho :\gd \rightarrow \gl(n,\CC )$ be an admissible
representation. For any positive integer $m$ we denote by $\rho _m
:\gd \times  \gl(m,\CC )\rightarrow  \gl(mn,\CC )$ the representation
 given by $\rho _m (g^\vee ,r)=\rho  (g^\vee )\otimes r$. Assume that both  $\rho$ and  $\rho _m$
are {\it liftable}. For any  $\pi \in \Irr (G)$ we define rational
functions $\gam '_{\pi},\gam ''_{\pi}$ on $\Irr \gl(m,F)$ by
$\gam '_{\pi}(\sigma)\doteq \gam_{\rho _m}(\pi,\sigma),\gam
''_{\pi}(\sigma)\doteq \gam _{n,m}(l_{\rho} (\pi),\sigma)$

\prop{} For any $\pi \in \Irr (G)$ the two rational functions
 $\gam '_{\pi}$ and $\gam ''_{\pi}$ on $\Irr \gl(m,F)$ coincide.
\eprop

This Proposition follows  from the main result of \cite{Har}.
\cor{kaz-gelf}
For any  $\rho :\gd\rightarrow \gl(n,\CC)$ such that the
representation $\rho
_{n-1}$ is  liftable,  the lifting
$l_{\rho}$   is determined by the knowledge
of the function $\gam_{\rho _{n-1}}(\pi,\star)$ on $\Irr
\gl(n-1,F)$. Moreover it is sufficient to know the  function $\gam_{\rho
_{n-1}}(\pi,\star)$ only  up to a multiplication by a
constant.
\ecor


\sec{torus}{The case of split tori}
In this section we assume that $F$is non-archimedian.

\ssec{}{Lifting for split tori}Almost the only case when $l_\rho$ is known
for every representation $\rho$ is the case when $\bfG$ is equal
to a split torus $\bfT$. In this section we show how to carry out our
program in this easy case.

We start with an explicit  description of the lifting $l_\rho$.

We can write $\rho$ in the form

\eq{}
\rho=\bigoplus\limits_{i=1}^n \lam_i
\end{equation}
where $\lam_i$ are characters of $\td$. By the definition of  $\td$ we
can  regard  $\lam_i$ as  cocharacters $\lam_i:\GG_m\to\bfT$ of $\bf T$.
Let $\Lam =X_{*}(\bfT)$ be the group of cocharacters of $\bfT$. We
define $\text{supp} (\rho)\subset \Lam$ as the union of $\lam _i ,1\leq i\leq n$.

Let $\theta$ be a character of $T$ and $\chi_i=\theta\circ\lam_i$
be the corresponding characters of $F^*$. Then for every
$i=1,...,n$ there exists $z_i\in\mathbb R _{>0}$ such that
$|\chi_i(t)|=z_i^{v_F(t)}$ for every $t\in F^*$.

Let us assume that the ordering $\lam_1,...,\lam_n$ is chosen
in such a way that
\eq{}
z_1\leq z_2\leq ...\leq z_n
\end{equation}

Let $B_n$ be the subgroup of upper triangular matrices
in $\gl(n,F)$. We can define the character $\chi:B_n\to\CC^*$
by setting
\eq{}
\chi(b)=\prod\limits_{i=1}^n \chi_i(b_{ii})
\end{equation}
where $b_{ii}$ are the diagonal entries of $b$.

Consider the unitarily induced representation
$i_{B_n}^{GL(n,F)}\chi$. It is well-known (cf. \cite{BZ2})
that this representation has a unique irreducible quotient.
We define this quotient to be $l_\rho(\theta)$.  As follows from the
theory of Eisenstein series  this definition satisfies the conditions
of \refco{lifting}.

\ssec{gam-torus}{$\gam$-functions for split tori}Let $\bfT$ be as above,
\eq{}
\rho=\bigoplus\limits_{i=1}^n \lam_i
\end{equation}
 an admissible representation of $\td$.We will give now an explicit
description of the distribution $\Phi_{\rho,T}$.
 Let  $\td_\rho=\GG_m^n$ be the
corresponding torus in $\gln$. Thus we get a homomorphism
$\bfp_\rho^{\vee}:\td\to\td_\rho$.  Let $\bfT_\rho$
be the split torus over $F$ dual to $\td_\rho$
(thus $T_\rho=(F^*)^n$). We then get an
$F$-rational map $\bfp_\rho:\bfT_\rho\to\bfT$.

Consider the
top degree differential form $\ome_\rho :=dt_1...dt_n$ on  $\bfT_\rho$.
and the function
$\bff_\rho:\bfT_\rho\to\AA^1$ where
$\bff_\rho(t_1,...,t_n)=t_1+...+t_n$.

For every $r\in \RR_{\geq 0}$ we denote by $T_\rho(r)$ the set
$\{ t\in T_\rho| |f_\rho(t)|\leq r\}$.

\prop{}
\begin{enumerate}
\item
For every $r\in \RR_{\geq 0}$ and for every open compact subset
$C\subset T$ the integral
\eq{}
\Phi_\rho^r(C)=\int\limits_{t\in p_\rho^{-1}(C)\cap T_\rho(r)}
\psi(f_\rho(t))|\ome_\rho|
\end{equation}
is absolutely convergent.
\item
For every $C$ as above the limit
\eq{}
\lim_{r\to\infty}\Phi_\rho^r(C)
\end{equation}
exists. We denote this limit by $\Phi_\rho(C)$. We also denote by
$\Phi_\rho$ the corresponding distribution on $T$.
\item
For every character $\theta$ of $T$ there exists $s_0(\theta)\in\RR$
such that for every $s\in \CC$ such that $\RE(s)> s_0(\theta)$
the convolution $\Phi_\rho*\theta|\sig|^s$ is absolutely convergent
and
\eq{}
\Phi_\rho*\theta|\sig|^s=\gam_{\rho,T}\theta|\sig|^s
\end{equation}
\item
$\Phi_\rho$ is $\sig$-regular.

\end{enumerate}
\eprop

Our next goal is to describe explicitly the space  $\calS_\rho$ in the
case of a torus.

\ssec{}{The Fourier transform $\calF_\rho$}
We assume now that the representation $\rho$ is faithful.
This assumption will be kept until the end of Section 8.

Consider the map
$\calF_\rho: C_c^{\infty}(T)\to C^{\infty}(T)$ given by
\eq{}
\phi\mapsto |\sig|^{-1}\Phi_\rho*{^{\iota}\phi}
\end{equation}
where $^{\iota}\phi(x)=\phi(x^{-1})$.

Let also $L^2_\rho(T)=L^2(T,|\sig|d^*t)$.
The following lemma is straightforward.
\lem{}
$\calF_\rho$ extends to a unitary automorphism of 
$L^2_\rho(T)$.
\elem

We now set
\eq{}
\calS_\rho(T)=C^{\infty}_c(T)+\calF_\rho(C^{\infty}_c(T))\subset L^2_\rho(T)
\end{equation}
(the reader should compare this definition with the definition of a Schwartz
space $\calS(X)$ in \cite{BrK}).

\ssec{}{The function $\calC_\rho$}Let us give an example of a function
in $\calS_\rho(T)$. Let $T_0\subset T$ be the maximal compact subgroup of
$T$. Then we can identify the quotient $T/T_0$ with the lattice 
$\Lam$ of cocharacters of $\bfT$.

Let
\eq{}
\rho=\bigoplus\limits_{i=1}^n \lam_i
\end{equation}
be an admissible representation of $\td$. We denote the set of all $\lam_i$ above by 
$\text{supp}(\rho)$.

 Let us define a $T_0$-invariant
function $\calS_\rho$ on $T$ (i.e. a function on $T/T_0=\Lam$ by setting
\eq{}
\calC_\rho(\mu)=\#\{ (a_1,...,a_n)\in\ZZ_+^n|\ a_1\lam_1+...+a_n\lam_n=\mu\}
\end{equation}
\lem{}
$\calC_\rho\in\calS_\rho(T)$.
\elem
\ssec{tro}{The semigroup $\obT$}We now want to exhibit certain locality properties
of the space $\calS_\rho(T)$. For this we have to introduce first an additional
notation.
. Let $\Lam_\rho$ be the sub-semigroup
of $\Lam$  generated
by $\lam \in supp(\rho),\calO_\rho$ be the group algebra
of $\Lam_\rho$ over $F$ and $\obT :=\Spec \calO_\rho$.
It is easy to see that $\obT$ is an $F$-semi-group which contains $\bfT$
as an open dense subgroup.
\prop{tlocal}Let $\phi$ be a function on $T$. Then 
$\phi\in\calS_\rho(T)$ if and only if $\phi$ satisfies
the following conditions:
\begin{enumerate}
\item
The closure of $\text{supp}(\phi)$ in $\oT_\rho$ is compact.
\item
For every $x\in \oT_\rho$ there exists a neighbourhood
$U_x$ of $x$ in $\oT_\rho$ and a function $\phi'\in\calS_\rho(T)$
such that
\eq{}
\phi|_{U_x}=\phi'|_{U_x}
\end{equation}
\end{enumerate}
\eprop

\refp{tlocal} says that one can determine whether a function $\phi$ lies in 
$\calS_\rho(T)$ looking at its local behaviour around points of $\oT_\rho$.
We will discuss the notion of a local space of functions in \refs{local}
in more detail.

For every $x\in \obT$ we denote by $\bfZ_x$ the stabilizer of $x$
in $\bfT\x \bfT$. If $x\in \oT_\rho$  we define an
$Z_x$-module $\calS_{\rho,x}$ as the quotient
$\calS_\rho(G)/\calS^0_{\rho,x}(T)$ where
$\calS^0_{\rho,x}(T)$ the space of all functions from
$\calS_\rho(G)$, which vanish in some neighbourhood of
$x$.

\lem{}
There exists unique $Z_x$-invariant functional
$\eps_x:\calS_{\rho,x}\to \CC$ such that
$\eps_x( \calC_\rho )=1$.
\elem

Thus given $\phi\in\calS_\rho(T)$ we can produce a function
$\phi^\eps$ on $\oT_\rho$ setting
$\phi^\eps(x)=\eps_x(\phi)$. However, the function $\phi^\eps$
is not locally constant in general.
\ssec{}{Mellin transform}We now going to give yet another
(equivalent) definition of $\calS_\rho(T)$ using the
{\it Mellin transform} on $T$. 
Let $X(T)=\Hom(T,\CC^*)$. Then
$X(T)$ has the natural structure of an algebraic variety,
whose irreducible components are parametrized by
$\Hom(T_0,\CC^*)=\Hom(T_0,S^1)$. Choose
a Haar measure $d^*t$ on $T$. 

Let $\Lam$ be the lattice of cocharacters of $\bfT$. 
As before one can identify $\Lam$
with $T/T_0$ where  $T_0\subset T$ is the maximal
compact subgroup of $T$. We denote by
 $X^{\text{un}}(T)\subset X(T)$  the subset
of unitary characters of $T$.

We denote by $X_0(T)$ the component of $X(T)$, consisting of
characters whose restriction to $T_0$ is trivial.
One has the natural identification
$X_0(T)\simeq\td$ (where $\td$ denotes the dual torus to $\bfT$
over $\CC$).
For any $\lam \in supp(\rho)$ and choose a lift $\lam'$ of $\lam$ to
$T$.
Then we denote by $\bfp_\lam$ the regular function on
$X(T)$ defined by
$$
\bfp_\lam(\chi)=
\begin{cases}
\chi(\lam')\ \text{if $\chi\circ\lam|_{\calO^*}=1$}\\
0\qquad \text{otherwise}.
\end{cases}
$$

Clearly, this definition does not depend on the choice
of $\lam'$.

For $\phi\in\cic(T)$ we define   the Mellin transform
$M(\phi)$ of $\phi$ as a function  on $X(T)$ given by
\eq{}
M(\phi)(\chi)=\int\limits_T \phi(t)\chi(t)d^*t.
\end{equation}

The following lemma is standard.
\lem{}
 $M$ defines an isomorphism between $\cic(T)$ and
$\calO(X(T))$ which extends to an isomorphism between $L^2(T)$ and
$L^2(X^{\text{un}}(T))$.
\elem

Let $\rho$ be an admissible representation of $\td$ and
$\calF_{\psi,\rho}$ denote the operator defined above.
The following result is proved by a straightforward calculation.
\th{torus-main}
 A function $\phi\in L^2_\rho(T)$ lies in $\calS_\rho(T)$ if
and only if
$\underset{1\leq i\leq n}\prod(\bfp_{\lam_i}-1) M(\phi)$ is a regular
function on $X(T)$.
\eth

Let us now discuss the connection of the space  $\calS_\rho(T)$ with the
corresponding local $L$-functions.
 
\lem{}
\begin{enumerate}
\item For every $\phi\in\calS_{\rho}(T)$ and every character $\chi$ of
$T$ the integral
\eq{}
Z(\phi,\chi,s)=\int\limits_T \phi(t)\chi(g)|\sig(t)|^s d^*t
\end{equation}
is absolutely convergent for $\Re e(s)>>0$.
\item
$Z(\phi,\chi,s)$ has a meromorphic continuation to the whole of
$\CC$ and defines a rational function of $q^s$.
\item
The space of all $Z(\phi,\chi,s)$ as above is a finitely generated
non-zero fractional ideal of
the ring $\CC[q^s,q^{-s}]$. We shall denote this ideal by $J_{\chi}$.
Let also $L_\rho(s,\chi)$ to be the unique generator of $J_\chi$ of
the form $P(q^{-s})^{-1}$, where $P$ is a polynomial such that
$P(0)=1$.
\item
 $L_\rho(\chi ,s)=\prod\limits_{1\leq i\leq n}L(\chi\circ \lam_i ,s)$ where
 $L(\chi\circ \lam_i ,s)$ denotes the corresponding Tate's
L-function (cf. \cite{tate}).
\end{enumerate}
\elem



\sec{local}{Local spaces of functions}
We now want to generalize some constructions of the preceding Section
to the case of arbitrary reductive group $G$. Let us remind that we
assume
now that the representation $\rho$ is admissible and faithful.

In this section we assume that  we have constructed 
the corresponding distribution $\Phi_{\psi,\rho}$ which
satisfies the assumptions of \refs{gamma}.

\ssec{sat}{Saturations}
Let $\oX$ be a topological set and $X\subset\oX$ an open dense subset.
Given a space $L$ of functions $\oX$ we say that
$L$ is local if there exists a sheaf  $\calL$ on $\oX$ and  an embedding
$\calL\in i_*(C(X))$ such that $L=\Gam_c(\oX,\calL)$ where $C(X)$ is the
sheaf of continuous functions on $X$.

It is easy to see that when the space $X$ is totally disconnected
then locality of $L$ is equivalent to the following two conditions:

1) $C_c(X)\subset L\subset C(X)$

2) $L$ is closed under multiplication by elements of $C(\oX)$.

Given a subspace $V\subset C(X)$ we denote by $\calP_V$ the presheaf
on $\oX$ such that
$\Gam(U,\calP_V)$ is the subspace of $\phi\in C(U\cap X)$ consisting of
all $\phi\in C(U\cap X)$ such that there exists $v\in V$ for which
$v|_{U\cap X}=\phi$. Let $\calL_V$ be the associated sheaf. It is clear
that we have the natural embedding
$\calL_V\hookrightarrow i_*(C(X)) $.
We set $\oV=\Gam_c(\oX,\calL_V)$ and call $\oV$ the $saturation$
of $V$ with respect to $\oX$.
\ssec{}{$C^{\infty}$-version}
Let now $\oX$ be a closed subset of a $C^{\infty}$-manifold
$Y$. Let $d(\cdot,\cdot)$ be a metric on $\oX$ coming from a Riemannian
metric on $Y$.

Let $i:X\hookrightarrow \oX$
be an open dense embedding of a $C^{\infty}$-manifold $X$ into $\oX$.
Let also $\calL\subset i_* C^{\infty}(X)$
be a subsheaf. We denote by $\Gam_c^{as}(\oX,\calL)$ the
space of $C^{\infty}$-functions $\phi$ on $X$ such that
for any $x\in\oX$ and any $N>0$
there exists an open neighbourhood $U$ of $x$ and
a section $l\in \Gam(U,\calL)$ such that for any $y\in U$ one has
$$
\phi(y)-l(y)\leq C d(x,y)^N
$$
where $C$ is a constant (i.e. $C$ does not depend on $y$).

\medskip
\noindent
{\bf Example.} Let $X=\oX=(0,1)$ and let $\calL$ be the sheaf of polynomial
functions on $X$. Then $\Gam_c^{as}(X,\calL)=C^{\infty}_c(0,1)$.

\medskip
\noindent
Given a subset $V$ of $C^{\infty}(X)$ we can define its saturation
$\oV$ as $\Gam_c^{as}(\oX,\calL_V)$ where $\calL_V$ is defined as in
\refss{sat}.
\ssec{}{The Fourier transform $\calF_\rho$}
In
this subsection we would to define certain "twisted" analog of the Fourier
transform in the space of functions of $G$ attached to every $\rho$ as
above.
As before we  denote by $L^2_\rho(G)$ the space of $L^2$-functions on
$G$ with
respect to the measure $|\sig|^{l+1}|dg|$ where $l$ is the semi-simple
rank of $\bfG$ and $|dg|$ is a Haar measure on $G$.

For every $\phi\in C^{\infty}_c(G)$ we define a new function $\calF_\rho(\phi)$
\eq{}
\calF_\rho(\phi)=|\sig|^{-l}(\Phi_{\psi,\rho}*{^{\iota}\phi})
\end{equation}
where $^\iota\phi(g)=\phi(g^{-1})$.
In what follows we will assume the validity of the following conjecture.

\medskip

\noindent
\conj{unitary}
1. $\calF_\rho$ extends to a unitary operator on the space
$L^2_\rho(G)$ and $\calS_\rho$ is $\calF_\rho$-invariant.
\econj

\medskip

\noindent
{\bf Example.} Let $G=\gl(n,F)$ and $\rho$ be the standard representation of
$\gd=\gl(n)$. In this case $L^2_\rho(G)$ is the same as $L^2(M_n)$ with respect
to the additive measure on $M_n$ and $\calF_\rho$ coincides with the Fourier
transform on $M_n$ (where we identify $M_n$ with the dual vector space by
means of the form $\la A,B\ra=\tr(AB)$).

\ssec{}{The semigroup $\obg$}We would like to embed $\bfG$ as an
open subset in a larger affine variety $\obg$ in a $\bfG\x
\bfG$-equivariant
way.

The variety $\obg$ is an algebraic semigroup, containing $\bfG$
as an open dense subgroup. It is well-known (cf. \cite{Vin}) that in
order to define such a semigroup one needs to exhibit a  subcategory of
the category of finite-dimensional
$\bfG$-modules, closed under subquotients, extensions
and tensor products.

We will  define now the category of {\it $\rho$-positive}
representations of $\bfG$ which will satisfy the above properties.
Let $\lam:\td\to \GG_m$ be a character, which has non-zero multiplicity
in $\rho$. We can regard $\lam$ as a cocharacter of $\bfT$.
We say that a $\bfG$-module $\bfV$ is $\rho$-positive if for every
$\lam$ as above the representation $\rho\circ \lam$ of $\GG_m$ is
isomorphic to a direct sum of characters of the form
$t\mapsto t^i$ for $i\geq 0$. It is clear that the category of
$\rho$-positive representations satisfies the properties discussed
above and thus defines a semi-group $\obg$.
\ssec{}{The space $\calS_\rho(G)$}
Let 
\eq{}
V_\rho=C^{\infty}_c(G)+\calF_\rho(C^{\infty}_c(G))\subset L^2_\rho(G)
\end{equation}
Unlike in the case of the torus, one can show that the
space $V_\rho$ is almost never local (this is not so for example
when $\bfG=\gl(2)$ and $\rho$ is the standard representation of
$\gd\simeq\gl(2,\CC)$).
We define $\calS_\rho(G)$ to be the saturation of $V_\rho$ (the definition
makes sense both for archimedian and non-archimedian $F$).

\ssec{}{The function $\calC_\rho$}Let us now exhibit certain explicit element
in $\calS_\rho(T)$.
In what follows we choose a square root $q^{1/2}$ of $q$.
\sssec{}{The Satake transform}Let $K=\bfG(\calO)$ and let
$\calH_K=C^{\infty}_c(K\backslash G/K)$ be the corresponding Hecke
algebra. Recall that
there exists the natural isomorphism $S:\calH_K
\simeq \calO(\td)^W\simeq \calO(\gd)^{\gd}$, where $\calO(\td)^W$
is the algebra of $W$-invariant regular functions on $\td$. This
isomorphism
is characterized as follows. One can identify $\td$ with the group
of unramified characters of the torus $T$. For every
$\lam\in\td$ let $i(\lam)$ denotes the corresponding representation
of $G$ obtained by normalized induction of the character $\lam$
from $B$ to $G$. Let $v_\lam\in i(\lam)$ be the (unique up
to a constant) non-zero K-invariant vector in $i(\lam)$. Then
for every $\phi\in C^{\infty}_c(K\backslash G/K)$ one has
\eq{}
(\phi |dg|)\cdot v_\lam=S(\phi)(\lam)v_\lam
\end{equation}
\sssec{}{The function $\calC_\rho$}Let $\rho$ be as above. Let
us also denote by $2\del_\bfG^{\vee}$ the sum of positive coroots
of $\bfG$.
 Set
\eq{}
f_{i,\rho}(\lam)=\tr(\lam \cdot 2\del_\bfG^{\vee}(q^{-1/2}),\text{Sym}^i\rho)
\end{equation}
for every $i\in\ZZ_+$.

We define now
\eq{crho}
\calC_\rho=\sum\limits_{i=0}^{\infty}S^{-1}(f_{i,\rho})
\end{equation}
It is easy to see that the non-degeneracy assumptions on
$\rho$ imply that $\calC_\rho$ is a well-defined locally constant
function on $G$, which however does not have
compact support.

\smallskip
\noindent

\lem{}
$\calC_\rho\in \calS_\rho(G)$ and $\calF_\rho(\calC_\rho)=\calC_\rho$.
\elem

The lemma follows easily from \refco{compatibility} from \refs{gamma}. 

\conj{eps}

There exists unique up to a constant $Z_x$-invariant functional
$\eps_x:\calS_{\rho,x}\to \CC$.

\econj

\medskip
\noindent
{\bf Example.} Let $\bfG=\gln$ and take $\rho$ to be the standard
representation of $\gd\simeq \gln$. Then
$\obg=\bfM_n$ -- the semi-group of $n\x n$-matrices.
Also one has $\calS_\rho(G)=\calS(M_n)$ -- the space of
Schwartz-Bruhat functions on $M_n$. In this case $\calS_{\rho,x}$
is
one-dimensional for every $x$ and the functional $\eps_x$ is equal to
the
evaluation of a locally constant function at $x$.
\ssec{llocal}{Local $L$-functions}Here we give a definition of the
local $L$-factor for every admissible non-singular representation $\rho$.

Let $\pi$ be an irreducible representation of $\bfG$. We denote by
$M(\pi)\subset C^{\infty}(G)$ the space of all matrix coefficients
of $\pi$.

\conj{llocal}
\begin{enumerate}
\item
For every $\phi\in\calS_\rho(G)$ and every $m\in M(\pi)$ the integral
\eq{integral}
Z(\phi,m,s)=\int\limits_G \phi(g)m(g)|\sig(g)|^s |dg|
\end{equation}
is absolutely convergent for $\Re e(s)>>0$.
\item
$Z(\phi,m,s)$ has a meromorphic continuation to the whole of
$\CC$ and defines a rational function of $q^s$.
\item
The space of all $Z(\phi,m,s)$ as above is a finitely generated
non-zero fractional ideal of
the ring $\CC[q^s,q^{-s}]$. We shall denote this ideal by $J_{\pi}$.
\end{enumerate}
\econj

We now define $L_\rho(s,\pi)$ to be the unique generator of $J_\pi$ of
the form $P(q^{-s})^{-1}$, where $P$ is a polynomial such that
$P(0)=1$.

\sec{poisson}{The Poisson summation formula}

Let as before $K$ be a global field and $\calP(K)$ be the set
of places of $K$. We also denote by
$\AA_K$ the adele ring of $K$.

Let $\bfG$ be a split reductive algebraic group over $K$ and
 $\rho:\gd\to \gl(n,\CC)$ be as before. For every place
$\grp\in \calP(K)$ we can consider the Schwartz
space $\calS_{\rho.\grp}=\calS_\rho(\bfG(F_\grp))$ with the
distinguished function $\calC_{\rho,\grp}$ in it. We now define the
space $\calS_\rho(\gaf)$ as the restricted tensor product
of the spaces $\calS_{\rho,\grp}$ with respect to the
functions $\calC_{\rho,\grp}$.

Choose now a non-trivial character $\psi:\AA_K\to \CC^*$ such that
$\psi|_{K}$ is trivial. Then $\psi$ defines an additive
character $\psi_{\grp}$ of $K_\grp$ for every $\grp$.
The Fourier transform
\eq{}
\calF_{\rho,\psi}=\prod\limits_{\grp\in\calP(K)}\calF_{\rho,\grp}
\end{equation}
acts on the space $\calS_\rho(\gaf)$.

Every $\phi\in\calS_\rho(\gaf)$ gives rise to a function on
$\gaf$.

\conj{poisson}

a) There exists a $\bfG(K)$-invariant functional
$\eps:\calS_\rho(\gaf)\to \CC$ such that

1) $\eps(\phi)=\eps(\calF_\rho(\phi))$ for any $\phi\in\calS_\rho(\gaf)$.

2) Let $\phi=\prod \phi_\grp\in\calS_\rho(\gaf)$. Assume that there exists
a place $\grp_0\in\calP(K)$ such that $\phi_{\grp_0}\in\calH(G)$.
Then
\eq{}
\eps(\phi)= \sum\limits_{g\in \bfG(K)}\phi(g)
\end{equation}

b) The functional $\eps$ is supported on $\obg(K)$. In other words,
assume that we are given $\phi \in \calS_\rho(\gaf)$ such that
for every $g\in\obg(K)$ there exists a neigbourhood $U_g$ of $g$ such
that $\phi(x)=0$ for every $x\in U_g\cap \bfG(K)$. Then $\eps(\phi)=0$.
\econj

When $\bfG=\gln$ and $\rho$ is the standard representation,
the statement of \refco{poisson} is the standard Poisson summation
formula. In this case $\phi$ is a smooth function on the
space $\bfM_n(\AA_K)$ of adelic $n\x n$ matrices and one has
\eq{}
\eps(\phi)=\sum\limits_{g\in \bfM_n(K)}\phi(g)
\end{equation}

We do not know any explicit formula for $\eps$ in general.

Let now $\pi$ be an irreducible automorphic representation of $\gaf$,
$\pi_\grp$ -- its component at $\grp\in\calP(S)$. One can consider the
global L-function
\eq{}
L_\rho(\pi,s)=\prod\limits_{\grp\in\calP(S)}L_\rho(\pi_\grp,s)
\end{equation}

Using  the arguments of \cite{GJ} is easy to see that the validity of
the above conjectures
implies the meromorphic continuation and functional equation
of the corresponding automorphic L-functions.

\sec{algint}{``Algebraic'' integrals over local fields}
In the next two sections we assume that
the local field $F$ has characteristic $0$.

Let $\bfX$ be an algebraic
variety over $F$ of dimension $d$, $X=\bfX(F)$. Recall that
(cf. \cite{weil}) any $\ome\in\Gam(\bfX,\Ome^{\top})$ defined
over $F$ defines a measure $|\ome|$ on $X$.
\ssec{}{Locally integrable functions} We say that a function
$\phi:X\to\CC$ is {locally integrable} if for every
point $x\in X$ and any top-form
$\ome\in\Gam(\bfU,\Ome^d)$ defined in a compact neighbourhood ${\bf
U\subset X}$ of $x$ the integral
\eq{}
\int\limits_U |\phi(y)||\ome(y)|
\end{equation}
is convergent.
\th{vadik}(V.~Vologodsky) Let $\bff:\bfX\to \AA^1$ be a proper morphism and
let $\ome\in \Gam(\bfX, \Ome^d)$. Suppose that both are defined over
$F$. Then  there exist an
open compact subgroup $K$ of $\calO_F^*$ and a $K$-invariant function
$\alp(x)$ such that  distribution $f_!(|\ome|)$ on $F$ is equal to
$\alp(x)dx$ for $|x|>>0$.
\eth

\reft{vadik} is proved in the Appendix.

\cor{cvadik}Let $\psi:F\to \CC^*$ be a non-trivial character. Then
the integral
\eq{}
\int\limits_X \psi(f)|\ome|:=\int\limits_F f_!(|\ome|)\cdot \psi:=
\lim\limits_{N\to \infty}\int\limits_{|x|<N}f_!(|\ome|)\psi
\end{equation}
is convergent.
\ecor
\ssec{}{Algebraic-geometric distributions}Let as before $\bfX$ be an
algebraic variety over $F$.
\defe{alg-dist}By an algebraic-geometric distribution on $\bfX$ we
mean
a quadruple $\bPhi=(\bfY,\bff,\bfp,\ome)$ where

1) $\bfY$ is a smooth algebraic variety over $F$

2) $\bfp:\bfY\to\bfX$, $\bff:\bfY\to\AA^1$ are morphisms defined over
$F$

3) $\ome\in\Gam(\bfY,\Ome_\bfY^{\top})$
\edefe

By an isomorphism between two algebraic-geometric distributions
$\bPhi=(\bfY,\bff,\bfp,\ome)$ and $\bPhi'=(\bfY',\bff',\bfp',\ome')$  we
shall mean a {\it birational} isomorphism between
$\bfY$ and $\bfY'$ preserving all the structure.

Given an algebraic-geometric distribution $\bPhi$ we can try to define
a usual distribution $\Phi=\Phi_\psi$ on $X$ by writing
\eq{dist}
\Phi=p_!(\psi(f)|\ome|)=\lim\limits_{N\to\infty}p_!(\psi(f)|\ome|)|_{\{y\in
Y|\ |f(y)|\leq N\}}
\end{equation}
\refc{cvadik} guarantees that in a number of cases \refe{dist} makes sense.
For example one can show that the following result holds.
\prop{finite}Let an algebraic-geometric distribution $\bPhi$ as above be
given.
Assume that there exists an open dominant embedding $\bfY\hookrightarrow
{\widetilde \bfY}$ such that ${\widetilde \bfY}$ is smooth and such
that
$\bfp,\bff$ and $\ome$ extend to ${\widetilde \bfY}$. Assume
furthermore that the morphism
$\bfp\x\bff:{\widetilde \bfY}\to\bfX\x \AA^1$ is proper.
Then $\Phi$ given by \refe{dist} is a well-defined distribution on
$X$.
\eprop

We say that an algebraic-geometric distribution $\bPhi$ is
{\it finite} if it satisfies the conditions of \refp{finite}.
Also, we call $\bPhi$ {\it weakly finite} if it is finite over
the generic point of $\bfX$. We also say that $\bPhi$ is {\it
analytically
finite} if it is weakly finite and for any finite extension $E/F$ the
corresponding distribution $\Phi$ defined {\it a priori} on a dense
open subset of $\bfX(E)$ is locally $L^1$ on the whole of $\bfX(E)$.
In this case we call $\Phi$ the {\it materialization} of $\Phi$.
\lem{tor-finite}
Let $\bfX=\bfT$ be a split torus over $F$ and and let $\rho$ be an
admissible representation of $\td$. Let
$\bPhi_{\rho,\bfT}=(\bfT_\rho,\bff_\rho,\bfp_\rho,\ome_\rho)$ be given
by the formulas of \refss{gam-torus}. Then $\bPhi_{\rho,\bfT}$ is
analytically
finite and the distribution $\Phi _{\rho,\bfT}$ is the materialization
of  $\bPhi_{\rho,\bfT}$.
\elem
\ssec{reduction}{Reduction of algebraic-geometric distributions}
Let $\bPhi=(\bfY,\bff,\bfp,\ome)$ be an algebraic-geometric distribution and let
$\bfV$ be a vector-group defined over $F$.
We say that $\bfV$ acts on $\bPhi$ if we are given a free
action of $\bfV$ on $\bfY$ such that

1) $\bfp$ and $\ome$ are $\bfV$-invariant.

2) Let $\bfq:\bfY\to\bfZ=\bfY/\bfV$ be the quotient map.
Then for any $z\in\bfZ$ the restriction of $\bff$ to
$\bfq^{-1}(z)$ is an affine function.

Let $\bfV^*$ denote the dual vector space
to $\bfV$. We denote by $\bfs:\bfZ\to\bfV^*$ the map
sending every $z\in\bfZ$ to the linear part of $\bff|_{\bfq^{-1}(z)}$.
Set
\eq{}
\obY=\{ z\in\bfZ|\ \bfs(z)=0\}.
\end{equation}

We say that an action of $\bfV$ on $\bPhi$ is non-degenerate if
for any $\oy\in\obY$ the differential
$d\bfs: T_{\oy}\obY\to \bfV^*$ is onto. In this case the
variety
$\obY$ is smooth and for any $\oy\in\obY$ we get the natural
isomorphism $T_{\oy}\bfZ/T_{\oy}\obY\simeq \bfV^*$. Hence for any
$\oy\in\obY$ and any $y\in\bfq^{-1}(\oy)$ we have the natural isomorphisms
\begin{align*}
\Lam^{\top}(T_{y}\bfY)\simeq
\Lam^{\top}(T_{\oy}\bfZ)\ten \Lam^{\top}(\bfV)\simeq
\Lam^{\top}(T_{\oy}{\obY})\ten \Lam^{\top}(T_{\oy}\bfZ/T_{\oy}{\obY})\ten\Lam^{\top}(\bfV)\simeq\\
\Lam^{\top}(T_{\oy}{\obY})\ten \Lam^{\top}(\bfV^*)\ten \Lam^{\top}(\bfV)\simeq
\Lam^{\top}(T_{\oy}\obY).
\end{align*}
Therefore, for any $\oy \in\obY, y\in\bfq^{-1}(\oy)$,
the restriction of $\ome$ to $T_y{\bfY}$ defines an element
$\oome_{\oy}\in\Lam(T_{\oy}(\bfY))^*$. Since $\ome$ is
$\bfV$-invariant, it follows that
$\oome_{\oy}$ does not depend on the choice of
$y\in\bfq^{-1}(\oy)$. It is easy to see that there exists a top-degree
form $\oome$ on $\obY$ such that
for any $\oy\in\obY$ one has $(\oome)_{\oy}=\oome_{\oy}$.

By the definition, the restriction of $\bff$ to
$\bfq^{-1}(\bfY)$ is $\bfV$-invariant. So,
$\bff|_{\bfq^{-1}\obY}=\obf\circ\bfq$ for some function $\obf$ on $\obY$.
Also, since $\bfp$ is $\bfV$-invariant we have
$\bfp|_{\bfq^{-1}(\obY)}=\obp\circ\bfq$ for some morphism
$\bfp:\obY\to\bfX$.

We call the algebraic-geometric distribution
$\obPhi=(\obY,\obp,\obf,\oome)$ the {\it reduction} of $\bPhi$.

Let us see the effect of this reduction procedure on materializations.
Assume that $\obPhi$ is analytically finite and denote by
$\oPhi$ its materialization. For simplicity we shall assume that
$F$ is non-archimedian. Let $V_n$ be sequence of open
compact subgroups of $V$, such that
\eq{}
V=\bigcup\limits_{n=1}^{\infty}V_n
\end{equation}
Let also $Y_n$ be an increasing covering of $Y$ by $V_n$-invariant open
compact subsets. Consider the sequence $\Phi_n$, $n=1,2,...$ of
distributions on $X$ where
\eq{}
\Phi_n(\phi)=\int\limits_{Y_n}\psi(f)p^*(\phi)|\ome|.
\end{equation}
\lem{reduction}The sequence $\Phi_n$ weakly converges to $\oPhi$.
\elem

When $\dim\bfV=1$ this is proved in \cite{GeKa}($\S 1$). The general case
is treated similarly.
\ssec{weq}{$W$-equivariant distributions}Let now $\bfG$ be a reductive
algebraic group over $F$ and let $\bfT$ be its Cartan group. The Weyl
group $W$ acts naturally on $\bfT$ and we have the natural
$F$-rational morphism $\bfs:\bfG\to\bfT/W$. Let $\bfG_r$ denote the
set
of regular elements in $\bfG$. Then the restriction of $\bfs$ to
$\bfG_r$ is
a smooth map. In the rest of this section $\bfs$ will always denote
the
restriction of $\bfs$ to $\bfG_r$.

Let now $\bPhi_\bfT=(\bfY_\bfT,\bff_\bfT,\bfp_\bfT,\ome_\bfT)$
be an algebraic-geometric distribution on $\bfT$. We assume that
$\bfp_\bfT$ is a dominant map.
By a {\it $W$-equivariant structure} on $\bPhi$ we shall mean
a birational action of $W$ on $\bfY$ such that

1) $\bff_\bfT$ is invariant under this action.

2) $\bfp_\bfT$ intertwines the action of $W$ on $\bfY_\bfT$ and on $\bfT$.

3) One has $w^*(\ome_\bfT)=(-1)^{l(w)}\ome_\bfT$, where $l:W\to\ZZ$ is the
length
function.

Since $W$ is a finite group we can always find an open dense subset of
$\bfY_\bfT$ on which the action $W$ is biregular. We can also assume that
this
subset is chosen in such a way that it
is contained in the set of regular semi-simple elements. Note that birational
modifications of $\bfY_\bfT$ (such as passing to an open subset) do not change
the distribution $\Phi$ given by \refe{dist}.
We also assume in the sequel
that
$\bfp_\bfT$ is a smooth map.

Assume now that we are given a $W$-equivariant algebraic-geometric distribution
$\bPhi_\bfT$ as above. Then we may construct an algebraic-geometric
distribution
$\bPhi=(\bfY_\bfG,\bff_\bfG,\bfp_\bfG,\ome_\bfG)$ over $\bfG_r$ in the
following way. We may assume that $W$ acts biregularly on $\bfY_\bfT$.
Then we define
\eq{}
\bfY_\bfG=\bfG_r\underset{\bfT/W}\times\bfY_\bfT/W
\end{equation}
We let $\bfp_\bfG$ be the projection on the first multiple.
Condition 1 above implies that $\bff_\bfT$ descends to a function
on $\bfY_\bfT/W$ and hence defines a function $\bff_\bfG$ on
$\bfY_\bfG$.

The differential form $\ome_\bfG$ on $\bfY_\bfG$ is defined as
follows.
Let $y\in\bfY_\bfG$ and let $g=\bfp_\bfG(y)$. Choose a Borel
subgroup $\bfB$ of $\bfG$ containing $g$ and let $\bfU$ be its
unipotent radical. Then we have the canonical isomorphism
$\bfB/\bfU\simeq \bfT$ and we set $t=g\text{mod}\bfU\in\bfT$. Let
$y=(g,y'')$ where $y''\in\bfY_\bfT/W$. Since $g$ is a regular and
semi-simple  there exists unique $y'\in\bfp_\bfT ^{-1}(t)$ such that
$y'\text{mod}W=y''$.

Let $\bfO_g$ denote the orbit of $g$ under the adjoint action.
 Then we have the short exact sequence
\eq{}
0\to T_{y'}\bfY_\bfT\to T_y\bfY_\bfG \to T_g\bfO_g\to 0
\end{equation}
where $T_y\bfY_\bfG,T_{y'}\bfY_\bfT,T_g\bfO_g$ are the corresponding
tangent spaces.

The choice of $\bfB$ induces a symplectic form $\alp_\bfB$ on
$T_g\bfO_g$. The tensor product
$\alp_\bfB^{\dim\bfO_g/2}\ten \ome_\bfT|_{T_{y'}\bfY_\bfT}$ defines
an element in $\Lam^{\dim\bfY_\bfG}( T_y\bfY_\bfG)^*$. It is easy to
see that condition 3 above implies that this element does not depend
on the choice of $\bfB$ and we define it to be
$\ome_{\bfG}|_{T_y\bfY_\bfG)}$.
\conj{MAIN}
\begin{enumerate}
\item
Let $\bfG$ be a split reductive group over $F$ and let
$\rho :\gd\to\gl(n,\CC)$ be an admissible representation of $\gd$. Let
$\bPhi_{\rho,\bfT}$ be the algebraic-geometric distribution as in
\refl{tor-finite}. Then there exists a $W$-equivariant structure on
$\bPhi_{\rho,\bfT}$ such that the corresponding distribution
$\bPhi_{\rho,\bfG}$ is analytically finite and the underlying
distribution on $G$ is equal to $\Phi_{\rho,G}$ introduced in
\refs{gamma}.

Let $\bfP\subset \bfG$ be a parabolic subgroup and let $\bfL$ be
 the
 corresponding Levi factor.  The representation $\rho$ defines
canonically a representation $\rho_\bfL$ of $\bfL^{\vee}$, with the
same restriction to $\td$. Note that the Weyl group $W_\bfL$ is
naturally
a subgroup in $W$.
\item The restriction of the $W$-action on $\bPhi_{\rho,\bfT}$ to
$W_\bfL$ is equal to the corresponding $W_\bfL$-action on
$\bPhi_{\rho_\bfL,\bfT}$.
\end{enumerate}
\econj

\noindent
{\it Remark.} We shall see some examples of the above $W$-equivariant
structure in some examples. Surprisingly, the corresponding $W$-action
on $\bfT_\rho$ turns out to be rather complicated and it is related
to the theory of geometric crystals (cf. \cite{BerK}).



\sec{action}{The case of $G(m,n)$ and its applications}
Let $\bfG$, $\bfT$ and $\rho$ be as above.
In \refs{torus} we have constructed a finite algebraic-geometric
distribution
$\bPhi_{\rho,\bfT}=(\bfT_\rho,\bff_\rho,\bfp_\rho,\ome_\rho)$
on $\bfT$. It is explained in
\refs{algint} that in order to get from it a finite
algebraic-geometric
distribution $\Phi_{\rho,\bfG}$ one has to define a birational
action of the Weyl group $W$ on $\bfT_\rho$ satisfying certain
conditions. We do not know how to construct this action in general.
In this section we give an explicit construction of this
action for the group $\bfG=\bfG(n,m)$ (cf. \refss{basex}).

\ssec{two-n}{The case $m=2$}

We now want to define the  action of $W$  on $\bfT_\rho$  in the case when
$m=2$. The  Weyl group
$W$ of $\bfG=\bfG(2,n)$ is isomorphic to $\ZZ_2\times S_n$.
In this subsection we are going to give an explicit formula
for the action of the first factor by using part 2 of \refco{MAIN}.

Let $\bfT$ be as before a maximal torus of $\bfG$. It can
be identified with the variety of
all collections $(a_1,a_2,b_1,...,b_n)$, where $a_i,b_j\in\GG_m$
and $a_1a_2=b_1...b_n$.

The torus $\bfT_\rho$ consists of all matrices $t=(t_{ij})$
where $i=1,2$ and $j=1,...,n$. The map $\bfp_\rho$ is given by
\eq{}
\bfp_\rho((t_{ij}))=(t_{11}t_{12}...t_{1n}, t_{21}t_{22}...t_{2n},
t_{11}t_{21},t_{12}t_{22},...,t_{1n}t_{2n})
\end{equation}
The function $\bff_\rho$ sends the matrix $(t_{ij})$ to
$\sum_{i,j}t_{ij}$ and the differential form $\ome_\rho$
is up to a sign equal to $\prod dt_{ij}$.

To define $\ZZ_2$-action on $\bPhi_{\rho,\bfT}$ we need to
define an involution $\tau$ on $\bfT_\rho$
that preserves $\bff_\rho$, sends $\ome_\rho$ to $-\ome_\rho$
and has the following two properties:

(1) For every $t\in\bfT_\rho$ and every $j=1,...,n$ one
has $t_{1j}t_{2j}=\tau(t)_{1j}\tau(t)_{2j}$.

(2) For every $t\in\bfT_\rho$
$t_{11}...t_{1n}=\tau(t)_{21}...\tau(t)_{2n}$ and
$t_{21}...t_{2n}=\tau(t)_{11}...\tau(t)_{1n}$.

It is easy to see that in the case when $n>2$ there exist many ways to
define an involution on $\bfT_\rho$ which satisfy the compatibility
conditions $(1)$ and $(2)$. We will use  the conjectural
properties of the distribution $\Phi _{\rho,\bfG(2,n)}$ to obtain the
formula for the involution $\tau$.

Let $\bfT_n=\GG_m^n$ and let $N:\bfT_n\to\GG_m$ be given
by $N((t_1,...,t_n))=t_1...t_n$.
Define
\eq{}
\bfG_n=\{ (g,t)\in\gl(2)\times \bfT_n|\ \det(g)=N(t)\}
\end{equation}
The group $\bfG_n$ can be naturally regarded as a Levi subgroup of
$\bfG(2,n)$. Hence the dual group $\gd_n$ is a Levi subgroup
of $\bfG(2,n)^{\vee}$ and we can restrict the representation $\rho$ to
it. This restriction can be explicitly described as follows.
The group $\gd_n$ is isomorphic to the quotient of
$\gl(2,\CC)\x (\CC^*)^n$ by the subgroup consisting of elements
of the form $(x\cdot\Id,x^{-1},x^{-1},...,x^{-1})$.
For every $i=1,...,n$ we can define a two-dimensional
representation $\rho_i$ of $\gd_n$, sending
$(g,z_1,...,z_n)$ to $gz_i$. Then it is easy to see that
\eq{sumrho}
\rho=\bigoplus\limits_{i=1}^n \rho_i.
\end{equation}

The Cartan group
of $\bfG_n$ is naturally isomorphic to the Cartan group
$\bfT$ of $\bfG(2,n)$ considered in \refss{two-n}.
The Weyl group of $\bfG_n$ is equal to the first factor  of the  Weyl
group $W$   of $\bfG(2,n)$.

As follows from
\refe{dir-sum} and \refe{sumrho} the distribution $\Phi_{\rho,G_n}$ is
equal to the convolution

\eq{}
\Phi_{\rho,G_n}(g,b_1,...,b_n)=(\Phi_{b_1}*...*\Phi_{b_n})|db_1...db_n|
\end{equation}
where  distributions
$\Phi_b ,b\in  F^*$ on $\gl(m,F)$ are defined in \refss{variant}.
In other words, $\Phi_{\rho,G_n}$ can be thought of as a materialization of
the algebraic-geometric distribution $\tbphi=(\tbY,\tbf,\tbp,\tilome)$,
where:

$\bullet$ $\tbY=\gl(2)^n$

$\bullet$ $\tbf(g_1,...,g_n)=\tr(g_1)+...+\tr(g_n)$

$\bullet$ $\tbp(g_1,...,g_n)=(g_1...g_n,\det(g_1),...,\det(g_n))$

$\bullet$ $\tilome=\det(g_1...g_n) dg_1...dg_n$

where  $dg$ denotes a translation invariant top degree
differential form on $\gl(2)$.

More precisely, for every $N>0$ define an open subset
$\tilY(N)$ of $\tilY$ by
\eq{}
\tilY(N)=\{ (g_1,...,g_n)\in \gl(2)^n|\ |(g_i)_{\alp\beta}|\leq N\}
\end{equation}
Let $\tilmu_N$ denote the restriction of the distribution
$\psi(\tilf)|\tilome|$ to $\tilY(N)$.
\lem{}
\eq{}
\Phi_{\rho,G_n}=\lim\limits_{N\to\infty} p_!(\tilmu_N).
\end{equation}
\elem

Fix a maximal unipotent subgroup $\bfU_+$ in $\gl(2)$. Then
$\bfU_+\simeq \GG_a$. For every $g\in\gl(2)$ the function
$u\mapsto \tr(gu)-\tr(g)$ is linear in $u$. We set
$\tr(gu)-\tr(g)=g_-u$ (thus $g\mapsto g_-$ is a function on $\gl(2)$).

Let $\bfV_+ =\bfU_+^{n-1}$. Define an action of $\bfV_+$ on $\tbY$ by
setting
\eq{}
(u_1,...,u_{n-1})(g_1,...,g_n)=
(g_1u_1^{-1},u_1g_2u_2^{-1},...,u_{n-1}g_n)
\end{equation}
This action preserves $\tbp$ and $\tilome$ and for
any $(g_1,...,g_n)\in\tbY$ we have
\eq{}
\tbf(ug)-\tbf(g)=\sum\limits_{i=1}^{n-1} u_i((g_{i+1})_- -(g_i)_-).
\end{equation}

Thus we obtain an action of $ \bfV_+$ on
$\tbphi=(\tbY,\tbf,\tbp,\tilome)$. Let $\obPhi _{U_+}=(\obY,
\obp ,\obf, \oome)$ be the $\bfU_+$-reduction of $\tbphi$ (cf. \refss{reduction}).
As follows from \refl{reduction} the  materialization ${\bf \bar \Phi _{U_+}}$ is equal
to $\Phi_\rho$.

Let $\mathcal B$ be the variety of Borel subgroups of $\gl(2)$. For
any $g\in \gl(2)$ we denote by $\mathcal B ^g \subset \mathcal B$
the subvariety  of Borel subgroups containing $g$.

\lem{}For any regular $g\in \gl(2)$ the centralizer of $g$ in
 ${\bf PGL}(2)$ acts simply transitively on the variety
$\mathcal B -\mathcal B ^g$.
\elem

\cor{}For any two unipotent subgroups $ {\bf U_+,U_-} \subset \gl(2)$
we have a canonical isomorphism $\obPhi _{\bfU_+}\simeq
\obPhi _{\bfU_-} $. Therefore, we can write $\obPhi$ instead of
 $\obPhi _{\bfU_+}$.
\ecor

Let
\eq{}
{\widetilde \bfG_n}=\{ (\bfB,g,t)\in\calB\x \bfG_n|\ g\in\bfB\}
\end{equation}
We have the natural maps $\bfr:\tbG_n\to \bfG_n$ (sending $(\bfB,g,t)$ to
$g$) and $\bfm:\tbG_n\to\bfT$ (sending $(\bfB,g,t)$ to $g\text{mod}\bfU_\bfB$
where $\bfU_\bfB$ denotes the unipotent radical of $\bfB$).
Then $r$ is a ramified double covering and for any $(g,t)\in
{\bfG_n}$ one has
$r^{-1}(g,t)= \mathcal B ^g$. We denote by $\tbphi
=(\tbY,\tbp ,\tbf ,\tilome )$ the lift of
${\bf \Phi}$ to  $\tbG_n$. By definition $\tbY=\obY
\underset {\bfG_n}\times \tbG_n$.

Since the
generic fiber of $\bfm$ carries a canonical top form (see \refss{weq}) we can define
an algebraic-geometric distribution $\bfm^{*}(\bPhi _{\rho , \bfT})$ on $\tbG_n$. We write
$\bfm^*(\bPhi _{\rho , \bfT})=(\bfY'_\bfT ,\bfp''_\bfT ,\bff''_\bfT ,\ome ''_\bfT )$, 
where
${\bf Y'_T =Y_{\rho , T}}\underset {\bfT}\times \tbG_n$. We will
construct now an isomorphism $\mu : \tbphi \to  m^{*}(\bPhi
_{\rho , \bfT})$.

For any $\bfB\in\calB$ such that $\bfU_+ \not\subset \bfB$ the
multiplication map gives rise to birational isomorphisms
${\bf B\times  U_+}\to \gl(2)$ and ${\bf U_+\times  B}\to \gl(2)$.
Therefore the subset $\tbY _+
:=\bfB^{n-1}\times  \gl(2) \subset \tbY$ is a section of the action of
$\bfV_+$ on  ${\tbY}$. Fix any element $(g,t;\bfB)\in  \tbG_n$ such
that $g$ is regular
and semisimple.
Since by definition $g\in \bfB$ we
see that $\tbY _+ \bigcap \tbp ^{-1}(g,t;\bfB)\subset\bfB^n$. So we
can describe the variety $\tbp ^{-1}(g,t;\bfB)$ as the preimage ${\bf s_B}
^{-1}(0)$  where  ${\bf s_B}$  is the restriction of the map $\bfs :
\tbY \to \bfV^*$ to  $\bfB^n \bigcap \bfp^{-1}(g,t)$ where $\bfs$
is as  in \refss{reduction}.
In other words we can identify the variety ${\tbY}$ with
${\bf s_B}^{-1}(0)$.

Let ${\bf U \subset B }$ be the unipotent radical of $\bfB$ and
let ${\bfT}_g\subset \bfB$
be the centralizer of $g$. Then   ${\bf B}=\bfT_g\bfU$ and  ${\bfT}_g$ is
canonically isomorphic to $\GG _m ^2$. Therefore  ${\bfT}_g ^n$ is
canonically isomorphic to $\bfT_\rho$.

\lem{}
a)The natural projection $\bfB\to  \bfT_g$ defines an
embedding $\mu: \bfs_\bfB^{-1}(0)\hookrightarrow  \bfT_\rho $ and the image
of $i$ is equal to ${\bf p_{\rho ,T}}^{-1}(g,t)$.

 In other words we have
constructed an isomorphism $\mu :\tbY \to \bfY'_\bfT$.

b)$\mu$ defines an isomorphism of algebraic-geometric distributions.
\elem

The natural involution $\theta$ of  $\tbG_n$ over  $\bfG_n$
induces an involution $\tiltet$ of  $\tbphi$ and we
can define  an involution  $\theta'$ of $\bfm^{*}(\bPhi _{\rho ,\bfT})$
by $\theta '=\mu \circ  \tiltet \circ \mu ^{-1}$. It is easy to
see that the involution  $\theta '$ is $\gl(2)$-invariant and
therefore it is  induced by an involution $\tau$ on  $\bPhi _{\rho
,\bfT}$.

Let us give an explicit formula for $\tau$.
 For every $k=1,...,n$ we define the function $\Del_k$ on
$\bfT_\rho$ by
\eq{}
\Del_k(t)=t_{11}...t_{1(k-1)}+t_{11}...t_{1(k-2)}t_{2k}+...
+t_{22}...t_{2k}
\end{equation}
(by definition $\Del_1(t)=1$).

We also define a rational function $\eta$ on $\bfT_\rho$
setting
\eq{}
\eta(t)=\frac{t_{11}...t_{1n}-t_{21}...t_{2n}}
{\Del_n(t)}.
\end{equation}

\lem{}The involution $\tau$ satisfies

\eq{tau-form}
\tau(t)_{21}\tau(t)_{22}...\tau(t)_{2k}=
t_{21}...t_{2k}+\Del_k(t)\eta(t)
\end{equation}
\elem
It is clear that $\tau$ is uniquely determined by
\refe{tau-form} and by the conditions (1) and (2) above.

Assume for example that $n=2$. Then one can compute the above action
explicitly. Namely, in this case we have
\eq{}
\bfT=\{(a_1,b_1,a_2,b_2)\in\GG_m^4|\ a_1b_1=a_2 b_2\}
\end{equation}
Let $w\in W$ be the involution which interchanges $a_1$ and $b_1$.
Let $\tau=\tau_w$. Then an explicit calculation shows that
\eq{two-by-two}
\tau:
\begin{pmatrix}
t_{11}\  t_{12}\\
t_{21}\  t_{22}
\end{pmatrix}
\mapsto
\begin{pmatrix}
t_{21}\frac{t_{11}+t_{22}}{t_{12}+t_{21}}\
t_{22}\frac{t_{12}+t_{21}}{t_{11}+t_{22}}
\\
t_{11}\frac{t_{12}+t_{21}}{t_{11}+t_{22}}\
t_{12}\frac{t_{11}+t_{22}}{t_{12}+t_{21}}
\end{pmatrix}
\end{equation}
Moreover, it is easy to see that
in this case $\tau$ given by \refe{two-by-two} is the
unique birational involution of $\bfT_\rho$ which satisfies our
requirements.
\ssec{gmn}{The case of $\bfG(n,m)$}Let us now consider the
group $\bfG(n,m)$ for arbitrary $n$ and $m$.
In this case one can identify $\bfT_\rho$ with the
variety $(t_{ij},i=1,...,m,j=1,...,n)$ of $m\times n$ matrices
with non-zero entries. We wish to define a birational action
of the Weyl group $W=S_m\x S_n$ on $\bfT_\rho$. For every
$\alp=1,...,m-1$ we define a birational involution
$\tau^1_\alp$ on $\bfT_\rho$ in the following way:

1) All rows of $\tau^1_\alp(t=t_{ij})$ except for the
$\alp$ and $\alp+1$st are equal to the corresponding
rows of $t$.

2) The $\alp$- and $(\alp+1)$st rows of $\tau^1_\alp(t)$ are
obtained from those of $t$ by means of \refe{tau-form}.

Similarly, for every $\beta=1,...,n-1$ we define
$\tau^2_\beta(t)=\tau^1_\beta(t')$ where $t'$ is the
transposed matrix to $t$.
The following lemma is straightforward.
\lem{}
The involutions $\tau^1_\alp$, $\tau^2_\beta$ commute with
$\bfp_\rho$, preserve $\bff_\rho$ and map $\ome_\rho$
to $-\ome_\rho$.
\elem
The following theorem is proven in \cite{BerK}($\S$ 6.2).
\th{ber-kaz}
The involutions $\tau^1_\alp$ and $\tau^2_\beta$ define a birational
action of the group $S_m\times S_n$ on $\bfT_\rho$.
In particular, $\tau^1_\alp$ commutes with $\tau^2_\beta$ for
any $\alp$ and $\beta$.
\eth

Despite the fact the formulas \refe{tau-form} are quite
explicit we did not manage to prove \reft{ber-kaz} directly.
The proof, given in \cite{BerK} uses the machinery
of geometric crystals.
\ssec{}{$\gam$-functions}Let us still assume that $\bfG=\bfG(n,m)$.
Recall that we have the natural character $\sig:\bfG\to \GG_m$
(which sends the pair $(A,B)$ of matrices to $\det(A)=\det(B)$.
Using \reft{ber-kaz} we can define an algebraic-geometric
distribution $\bPhi_{\rho,\bfG}$ over the generic point of $\bfG$.
We would like to check that it gives the "correct" answer.

Let $\pi_m,\pi_n$ be generic representations of respectively
$GL(m,F)$ and $GL(n,F)$. Then in \cite{JPSS}
H.~Jacquet, I.~Piatetskii-Shapiro and J.~Shalika define a
rational function $\Gam(\pi_m,\pi_n,s)$ of one complex variable $s$
(which depends on the choice of a character $\psi$)
The following result can be deduced from \cite{BerK}.
\th{}
\begin{enumerate}
\item
The materialization $\Phi_{\rho,G}$ of $\bPhi_{\rho,\bfG}$ defines
a $\sig$-compact distribution on the whole of $G$. We denote
by $\gam_{\rho}$ the corresponding rational function on the
set $\Irr(G)$.
\item
Assume that $n\geq m$.
Let $\pi_m$, $\pi_n$ be as above and let $\pi$ be an irreducible
representation of $G$ such that
$\Hom_G(\pi,\pi_m\ten\pi_n)\neq 0$.
Then
\eq{}
\Gam(\pi_m,\pi_n,s)=\text{sign}(\pi_n)\gam_\rho(\pi|\cdot|^s)
\end{equation}
where $\text{sign}(\pi)$ is the value of the central character of
$\pi_n$ at the matrix $-\Id_n\in GL(n,F)$.
\end{enumerate}
\eth
\ssec{}{Lifting from non-split tori}Let us recall the notations of
\refss{lift-torus}: let $E/F$ be a separable extension of $F$ of degree $n$.
Let $\alp_E:\text{Gal}(\oF/F)\to S_n$ be the corresponding
homomorphism.
We denote by  $\bfN:\bfT_E\to\bfG_m$  the morphism of algebraic
groups coming from the norm map
$N:E^*\to F^*$.

Let $\bfT_E=\Res_{E/F}\GG_{m,E}$ where $\Res_{E/F}$ denotes the
functor of restriction of scalars. For every $m>0$ we set
\eq{}
\bfG_{E,m}=\{ (t,g)\in \bfT_E\x \gl(m)|\ \bfN(t)=\det(g)\}
\end{equation}

As is explained in \refss{kaz-def}, in order to define the lifting
$l_E(\theta)\in\Irr(\gl(n,F))$ of a
character $\theta:E^*\to \CC^*$, it is enough to
define the $\gam$-function $\gam(l_E(\theta),\pi_m)$ as a function on
$\Irr (\gl(m,F))$. Moreover, it is enough to know this function only
up to a constant.

We are now going to give a
conjectural definition of $\gam(l_E(\theta),\pi_m)$ by
constructing a weakly finite algebraic geometric distribution
$\bPhi_{E,m}$ on $\bfG_{E,m}$ (in fact, we will define it
only up to a constant multiple). We conjecture that this distribution is
analytically finite and that the corresponding function
on $\Irr(E^*)\x \Irr(\gl(m,F))$ is equal to $\gam_{E,m}$.

Choose $d\in F^*$ such that its image in $F^*/(F^*)^2$ is equal to
the discriminant of $E$. Let us also choose its square
root $\sqrt d$ in $\oF$.

Let
\eq{}
\bfG_{E,m}'=\{(t_1,...,t_n,g)\in\GG_m^n\x \gl(m)|\ t_1...t_n=\det(g)\}
\end{equation}
Using \refss{gmn} we define an algebraic-geometric
distribution $\bPhi'=(\bfY',\bfp',\bff',\ome')$ on
$\GG_m^n\x \gl(m)$. Moreover, we have the natural
$S_n$-action on $\bfY'$, which is compatible with
$\bfp'$ and leaves $\bff'$ (resp. $\ome'$)
invariant (resp. skew-invariant).

Define a new $\text{Gal}(\oF/F)$-action on $\bfY'$ by
\eq{}
g^{\text{new}}(y)=\alp_E(g^{\text{old}}(y))
\end{equation}

This action defines a new $F$-rational structure on $\bfY'$.
We denote by $\bfY$ the corresponding $F$-variety.
It is clear that the map $\bfp'$ gives rise to an $F$-rational
morphism $\bfp:\bfY\to\bfG_{E,m}$. Also the function
$\bff'$ and the differential form $\sqrt d\ome'$ give rise
to $F$-rational function $\bff$ and differential form $\ome$
on $\bfY$. Thus we set
$\bPhi _{{\bf E},m}=(\bfY,\bfp,\bff,\ome )$ to be the required
algebraic-geometric distribution.  We denote by $\Phi _{{\bf E},m}$ the
materialization of $\bPhi _{{\bf E},m}$.

Clearly, for different choices of $\sqrt d$, the resulting
materializations
$\Phi _{{\bf E},m}$ will differ only by a multiplication by $c\in
\mathbb C ^{*}$. Hence, the above construction suffices in order
to determine the lifting uniquely.

 Since the local Langlands conjecture is
known (see \cite{Har}, \cite{HT} and \cite{LRS}) we can ask whether our
definition of lifting coincides with one which is implied by \cite{Har},
\cite{HT} and \cite{LRS}. More precisely let $\chi$ be a character of
the group $T={\bf T}(F)$. The local
class field theory associates to  $\chi$ a homomorphism $\theta _\chi :\text{Gal}
(\bar E /E)\to \CC ^{*}$. Let $\Theta _\chi :=
\Ind^{\text{Gal}(\oF/F)}_{\text{Gal}(\bar E /E)}\theta _\chi$.
Then  $\Theta _\chi$ is an
$n$-dimensional representation of the group $\text{Gal}(\oF /F)$.
 Since the local Langlands conjecture for  $\bfG=\gln$ is
known  one
associate with   $\Theta _\chi$ an
irreducible representation $\pi  _\chi \in \Irr (\gl(n,F))$. Therefore we
can consider the function $\Gam(\pi_m,\pi  _\chi ,s)$ on the set
$ \Irr\gl(m,F)$.
 One can ask whether there exists $c\in \mathbb C ^{*}$ such that for
any $\pi_m \in  \Irr\gl(m,F)$ we have  $\gam(\pi_m,\pi  _\chi ,0)=c \chi
\otimes \pi _m (\Phi _{{\bf E},m})$.

\sec{finite}{The case of finite fields}
In this section we shall a give an explicit conjectural
construction
of $\Phi_{\rho,G}$ for any $\bfG$ and for (almost) any
$\rho$ in the case when the field $F$ is finite.
It turns out that in this case the relevant tool
from algebraic geometry which allows to go from the
case of a torus to the case of arbitrary group is not
the language of algebraic-geometric distributions,
but that of $\ell$-adic perverse sheaves. This tool
allows to avoid constructing an action of $W$ on $\bfT_\rho$.
\ssec{}{Notations}
In this section we fix a finite field $F$ with $q$ elements,
a prime number $l$ different
from the characteristic of $F$, and  a non-trivial
character $\psi:F\to \qlb^*$. All representations discussed in the
note are over the field $\qlb$.  We denote by $\calL_\psi$ the
Artin-Schreier sheaf on $\AA^1_F$.

Let $\bfG$ be a reductive algebraic group over $F$ and
let $\bfT$ be its abstract Cartan group. The torus $\bfT$ comes
equipped with a canonical $F$-rational structure. In this section
we assume for simplicity that
$\bfG$ is split. Then this $F$-rational structure is split too.
We will denote by $\fr:\bfT\to\bfT$ the corresponding Frobenius
morphism. Let $W$ denote the Weyl group of $\bfG$. For
every $w\in W$ we set $\frw=w\circ\fr$. The morphism $\frw$ induces a new
$F$-rational structure on $\bfT$ and we will denote the
corresponding algebraic torus over $F$ by $\bfT_w$.
It is well-known that there exists an embedding
$\bfT_w\hookrightarrow \bfG$ and that in this way we get
a bijection between conjugacy classes of elements in $W$ and
conjugacy classes of $F$-rational maximal tori in $\bfG$.

For a character $\theta:T_w\to\qlb^*$ we denote by
$R_{\theta,w}$ the corresponding Deligne-Lusztig
representation.

For two $\ell$-adic complexes $A$ and $B$ on $\bfG$ we can define
the convolution complex
 $A\star B$ and by setting
\eq{}
A\star B=m_!(A\boxtimes B)[\dim \bfG]
\end{equation}
where $m:\bfG\x \bfG\to \bfG$ is the multiplication map.

Let $\bfX$ be an algebraic variety over $F$ and let $\calF$ be
a complex of $\ell$-adic sheaves on it. By
{\it a Weil structure} on $\calF$ we shall mean an isomorphism
$\xi:\fr^*\calF\to\calF$ where $\fr$ denotes the geometric
Frobenius morphism on $\bfX$. In this case we can define
a function $\Tr(\calF)$ on $X=\bfX(F)$ by setting
\eq{}
\Tr(\calF)(x)=\sum\limits_{i} (-1)^i \tr(\xi_x:H^i(\calF_x)\to
H^i(\calF_x))
\end{equation}
(here $H^i(\calF_x)$ denotes the $i$-th cohomology of the
fiber of $\calF$ at $x$).

Let us choose a square root $q^{1/2}$ of $q$. Then for
any Weil sheaf $\calF$ on $\bfX$ and any half
integer $n$ we can consider the Tate twist $\calF(n)$ of
$\calF$. By the definition one has
\eq{}
\Tr(\calF(n))=\Tr(\calF)q^{-n}
\end{equation}

\ssec{}{$\gam$-functions for $\gln$}
\sssec{}{The main result}
Let $(\pi,V)$ be an irreducible representation of $G=GL(n,F)$.
Choose a non-trivial additive character
$\psi:F\to \qlb^*$ as above and consider the operator
\eq{fing}
\sum\limits_{g\in G}\psi(\tr(g))\pi(g)(-1)^n q^{-n^2/2}\in \End_G V
\end{equation}
Since $\pi$ is irreducible, this operator takes the form
$\gam_{G,\psi}(\pi)\cdot \Id_V$ where $\gam_G(\pi)\in \qlb$
(we will omit the subscripts $G$ and $\psi$ when it does not lead to
a confusion).
The number $\gam(\pi)=\gam_{G,\psi}(\pi)$ is called the gamma-function
of the representation $\pi$. The purpose of this subsection is to
compute explicitly the gamma-functions of all irreducible
representations of $G$.

Let $W\simeq S_n$ denote the Weyl group of $\bfG$.

Fix $w\in W$. For a character $\theta:T_w\to \qlb^*$ we set
\eq{}
\gam_w(\theta)=(-1)^{n+l(w)} q^{-n/2}\sum\limits_{t\in T_w}\psi(\tr(t))\theta(t)\in\qlb
\end{equation}
\medskip
{\bf Example.} Assume that $w\in S_n$ is a cycle of length $n$.
Then $T_w\simeq E^*$ where $E$ is the (unique up to isomorphism)
extension of $F$ of degree $n$. In this case
$\gam_w(\theta)=\gam_E(\theta)$ for any character $\theta$ of
$E^*$, where by $\gam_E(\theta)$ we denote the $\gam$-function
defined as in \refe{fing} for the group $\gl(1,E)\simeq E^*$.
\th{finite-main}Assume that an irreducible representation
$(\pi,V)$ appears in $R_{\theta,w}$ for some
$w$ and $\theta$ as above. Then
\eq{}
\gam(\pi)=(-1)^{l(w)}\gam_w(\theta)
\end{equation}
where $R_+$ denotes the set of positive roots of $G$.
In particular, $\gam(\pi)=\gam(\pi')$ if $\pi$ and $\pi'$
appear in the same virtual representation $R_{\theta,w}$.
\eth
The rest of this subsection is occupied with the proof of
\reft{finite-main}.
\sssec{}{Character sheaves}Let $\bfG$ be an arbitrary
reductive algebraic group over $F$.
Let us recall Lusztig's definition of
(some of) the character sheaves. Let $\tbG$ denote the variety
of all pairs $(\bfB,g)$, where

$\bullet$\quad $\bfB$ is a Borel subgroup of $\bfG$

$\bullet$\quad $g\in \bfB$

One has natural maps $\alp: \tbG\to \bfT$ and $\pi:\ \tbG\to \bfG$
defined as follows. First of all, we set $\pi(\bfB,g)=g$.
Now, in order to define $\alp$, let us remind that for any
Borel subgroup $\bfB$ of $\bfG$ one has canonical identification
$\mu_\bfB:\ \bfB/\bfU_{\bfB}\wt{\rightarrow} \bfT$, where
$\bfU_{\bfB}$ denotes the unipotent
radical of $\bfB$
(in fact, this is how the abstract Cartan group $T$ is defined).
Now we set $\alp(\bfB,g)=\mu_\bfB(g)$.

Let $\calL$ be a tame local system on $\bfT$. We define
$\chl=\pi_!\alp^*(\calL)[\dim \bfG]$.
One knows (cf. \cite{lu-char}, \cite{laum})
that the sheaf $\chl$ is perverse.

Assume now, that for some $w\in W$ there exists an
isomorphism $\calL\simeq \frw^*(\calL)$ and let us fix it.
It was observed by
G.~Lusztig in \cite{lu-char} that fixing such an isomorphism endows
$\chl$ canonically with a Weil structure. Lusztig's definition of this
Weil structure was as follows.

Let $j:\bfG_{rs}\to \bfG$ denote the open embedding of the variety of
regular semisimple elements in $\bfG$ into $\bfG$.

\lem{}
\eq{}
\chl=j_{!*}(\chl|_{\bfG_{rs}})
\end{equation}

Here $j_{!*}$ denotes the Goresky-MacPherson (intermediate) extension (cf.
\cite{bgg}).
\elem

The lemma follows from the fact that the map $\pi$ is small in the
sense of Goresky and McPherson.

The lemma shows that it is enough to construct the Weil structure only
on the restriction of $\chl$ on $\bfG_{rs}$.
The latter now has a particularly
simple form. Namely, let $\tbG_{rs}$ denote the preimage of
$G_{rs}$ under $\pi$ and let $\pi_{rs}$ denote the restriction of
$\pi$ to $\tbG_{rs}$. Then it is easy to see that
$\pi_{rs}:\tbG_{rs}\to \bfG_{rs}$ is an unramified
Galois covering with Galois
group $W$. In particular, $W$ acts on $\tbG_{rs}$ and this action is
compatible with the action of $W$ on $T$ in the sense that the restriction
of $\alp$ on $\tbG_{rs}$ is $W$-equivariant.

Now, an isomorphism $\calL\simeq\frw^*(\calL)$ gives rise to an isomorphism
\eq{weil-glupost1}
\alp^*\calL\simeq(w\circ\fr)^*(\alp^*\calL)
\end{equation}
(here both $w$ and $\fr$
are considered on the variety $\tbG$). Since $\pi_{rs}$ is a Galois covering
with Galois group $W$, it follows that one has canonical identification
\eq{weil-glupost2}
\pi_{rs!}(\fr^*(\alp^*\calL))\simeq  \pi_{rs!}((w\circ\fr)^*(\alp^*\calL))
\end{equation}

Hence from \refe{weil-glupost1} and \refe{weil-glupost2} we get
the identifications
\eq{}
\fr^*\pi_{rs!}(\alp^*\calL)\simeq \pi_{rs!}(\fr^*(\alp^*\calL))\simeq
\pi_{rs!}((w\circ\fr)^*(\alp^*\calL))\simeq \pi_{rs!}(\alp^*\calL)
\end{equation}
which gives us a Weil structure on $\pi_{rs!}(\alp^*\calL)\simeq
\chl|_{\bfG_{rs}}$. Hence we have defined a canonical Weil structure on
$\chl$.

Let now $\theta:T_w\to \qlb^*$ be any character. It is well known
that one can associate to $\theta$ a one-dimensional local
system $\calL_{\theta}$ together with an isomorphism
$\calL_{\theta}\simeq \frw^*\calL_{\theta}$. This local system is constructed
as follows: consider the sheaf $(\frw)_*(\qlb)$. It is clear that this
sheaf has a natural fiberwise action of the group $T_w$. Thus we
set $\calL_{\theta}$ to be the direct summand of $(\frw)_*(\qlb)$ on which
$T_w$ acts by means of the character $\theta$.
The following result is due to G.~Lusztig.
\th{}
\eq{}
\Tr(\calK_{\calL_{\theta}})=(-1)^{\dim \bfG}\ch (R_{\theta,w})
\end{equation}
\eth
\sssec{}{}Now we assume again that $\bfG=\gln$.
Set $\phpsi=\tr^*\calL_\psi[n^2](\frac{n^2}{2})$.

In this case it is well-known that for every character
$\theta:T_w\to \qlb^*$ the vector space, spanned by
the characters of all irreducible constituents of
$R_{\theta,w}$ coincides with the vector space spanned
by the functions of the form $\Tr(A)$, where $A$ runs
over all direct summands of $\calK_{\calL_{\theta}}$.
Hence \reft{finite-main} follows from the following result.
\prop{convolution}
\begin{enumerate}
\item  Let $\tr_{\bfT}:\bfT\to\AA^1$ denote the restriction of the
 trace morphism
 from $\bfG$ to $\bfT$. Then for every $\calL$ as above one has
\eq{}
H^i_c(\tr_{\bfT}^*(\calL_\psi)\ten \calL)=
\begin{cases}
0\quad\text{if $i\neq n$},\\
\text{is one-dimensional if $i=n$}.
\end{cases}
\end{equation}
We set $H_{\calL}:=H^n_c(\tr_{\bfT}^*(\calL_\psi)\ten \calL)(\frac{n}{2})$.
\item
Let $\theta:T_w\to \qlb^*$ be a character (for some $w\in W$).
Then the natural isomorphism $\frw^*\calL_{\theta}\simeq \calL_{\theta}$
gives rise to a natural endomorphism of $H_{\calL_{\theta}}$ which
by abuse of language we will also denote by $\frw$. Then
\eq{}
\frw|_{H_{\calL_{\theta}}}=(-1)^{l(w)}\gam_w(\theta)\cdot\Id
\end{equation}
\item
One has
\eq{isom}
\phpsi\star\chl\simeq H_{\calL}\ten \chl
\end{equation}
Moreover, if $\calL$ is endowed with an isomorphism
$\frw^*\calL\simeq \calL$ then \refe{isom} is
$\frw$-equivariant.
\end{enumerate}
\eprop

The proof follows easily from standard properties of the
Fourier-Deligne transform. The details are left to the reader.
\ssec{}{The case of arbitrary group}
Let $\bfG$ be an arbitrary connected reductive algebraic group
over $F$. For simplicity we will assume that $\bfG$ is split.
Let $\gd$ denote the Langlands dual group of $\bfG$ (which is
a reductive algebraic group over $\qlb$).
Let $\rho:\gd\to \Aut(E)$ be a representation of
$\gd$, where $E$ is a vector space over $\qlb$ of dimension $n$.
We would like to associate to it a certain
function $\pi\mapsto \gam_\rho(\pi)$ on the set of isomorphism
classes of irreducible representations of $G$. This is done
as follows.

The group $\gd$ comes equipped with a canonical
maximal torus $\td$. Let us diagonalize $\rho|_{\td}$. Let
$\lam_1,...,\lam_n$ be the corresponding characters of $\td$
in $\rho$ (with multiplicities). Let
also $\bfM^{\vee}_\rho$ denote the minimal Levi subgroup
containing $\rho(\td)$ and such that the action of
$\bfM^{\vee}_\rho$ is multiplicity free. We denote by
$W_\rho'$ the group $\text{Norm}_{\Aut(E)}(\bfM^{\vee}_\rho)/
\bfM^{\vee}_\rho$. We also denote by $W_\rho$ the Weyl group
of $\Aut(E)$.
\lem{}$W_\rho'$ is naturally a subquotient of $W_\rho$.
\elem
\prf Standard.
\epr

Let $\td_\rho\simeq \GG_{m,\qlb}^n$
be the Cartan group of $\gln$. Thus we get a natural map
$\bfp_\rho^{\vee}:\td\to\td_\rho$ sending every $t$ to
$(\lam_1(t),...,\lam_n(t))$.

Let now $\bfT_\rho\simeq \GG_{m,F}$
denote the dual torus to $\td_\rho$ over $F$ and
let $\bfp_\rho:\bfT_\rho\to \bfT$ denote the map, which
is dual to $\bfp_\rho^{\vee}$. Explicitly one has
\eq{}
\bfp_\rho(x_1,...,x_n)=\lam_1(x_1)...\lam_n(x_n) .
\end{equation}

The representation $\rho$ defines the natural homomorphism
from $W$ to $W_\rho'$, which by abuse of the language we will
also denote by $\rho$.

Let now $\pi$ be an irreducible representation of $G$. Assume that
$\pi$ appears in some $R_{\theta,w}$ for some
$\theta:T_w\to \qlb^*$. Let $w'$ be any lift of
$\rho(w)$ to $W_\rho$. Then $\bfp$ induces an
$F$-rational map $\bfp_w:\bfT_{\rho,w'}\to\bfT_w$, hence
a homomorphism $p_w:T_{\rho,w'}\to T_w$. Define now
\eq{}
\gam_\rho(\pi):=\gam(\pi'),
\end{equation}
where $\pi'$ is any irreducible representation of
$G_\rho$ which appears in $R_{p_w^*(\theta),w'}$. By
\reft{finite-main} one has
\eq{}
\gam_\rho(\pi)=(-1)^{l(w')}\gam_{T_{w'}}(p_w^*(\theta)).
\end{equation}
\lem{}The definition of $\gam_\rho(\pi)$ does not depend on the
choice of $w'$.
\elem
Let now $\Phi_{\rho}$ denote the unique central function on
$G$ such that for every irreducible representation $(\pi,V)$ of
$G$ one has
\eq{}
\sum_{g\in G}\Phi_\rho(g)\pi(g)=\gam_{\rho}(\pi)\cdot\Id_V
\end{equation}
We would like to compute this function explicitly.
We will be able to solve only a slightly weaker problem.
Namely, we are going to construct explicitly an $\ell$-adic
perverse sheaf
$\phro$  on $\bfG$ such that
\eq{}
\Tr(\phro)=\tilPhi_\rho
\end{equation}
where $\tilPhi_\rho$ is a function which satisfies the following condition:
\eq{}
\sum\limits_{g\in G}\tilPhi_\rho(g)\pi(g)=\gam_\rho(\pi)\cdot\Id
\end{equation}
for every irreducible "generic" representation $\pi$ of $G$
(cf. \reft{phi-gamma} for the precise statement).

\sssec{}{The perverse sheaf $\phro$}In what follows we assume that the
zero weight does not appear in $\rho$.
Let $\tr_\rho:\bfT_\rho\to \AA^1$ be given by
\eq{}
\tr_\rho(x_1,...,x_n)=x_1+...+x_n.
\end{equation}
Consider the complex $A_\rho:=(\bfp_\rho)_!\tr_\rho^*\calL_\psi[n](\frac{n}{2})$ (on
$\bfT$). It is easy to see that this complex is perverse.
We would like to endow this complex with a $W$-equivariant structure.

Choose $w\in W$. We need to define an isomorphism
$\iota_w:w^*(A_\rho){\wt\to}A_\rho$. Let (as above)
$w'$ be any lift of $\rho(w)$ to $W_\rho$. Then one
has
\eq{bred}
\bfp_\rho(w'(t))=w(\bfp_\rho(t))
\end{equation}
The sheaf $\tr_\rho^*\calL_\psi$ is obviously $W_\rho$-equivariant.
This, together with \refe{bred} gives rise to an isomorphism
$\iota_w':w^*(A_\rho){\wt\to}A_\rho$.
We now define $\iota_w:=(-1)^{l(w')}\iota_w'$.
\prop{}The isomorphism $\iota_w$ does not depend on the choice of
$w'$. The assignment $w\mapsto\iota_w$ defines a $W$-equivariant
structure on the sheaf $A_\rho$.
\eprop

The above $W$-equivariant structure on $A_\rho$ gives rise to
the sheaf $B_\rho=(\bfq_! A_\rho)^W$ on the quotient $\bfT/W$
 where $\bfq:\bfT\to\bfT/W$ is the
canonical map.

Let $\bfG_r$ be the set of regular elements of $\bfG$ and let
$\bfi$ be its embedding to $\bfG$. Let $\bfs:\bfG_r\to\bfT/W$ be the morphism coming from the
identification $\bfT/W\simeq \calO(\bfG)^{\bfG}$.
We define
\eq{}
\phro:=i_{!*}\bfs^*(B_\rho)[\dim\bfG-\dim\bfT](\frac{\dim\bfG-\dim\bfT}{2})
\end{equation}
We also define
\eq{}
\tilPhi_\rho=\Tr(\phro)
\end{equation}
By the definition the complex $\phro$ on $\bfG$ is $\bfG$-equivariant
with respect to the adjoint action and therefore the function
$\tilPhi_\rho$ is central.

\sssec{}{The action of $\tilPhi_\rho$ in $R_{\theta,w}$}.
Recall that a local system $\calL$ on $\bfT$ is called
quasi-regular if for every coroot $\alp^{\vee}:\GG_m\to\bfT$
the local system $(\alp^{\vee})^*\calL$ is non-trivial.
We say that a character $\theta:T_w\to\qlb^*$ is
{\it quasi-regular} if the local system $\calL_\theta$
is quasi-regular.
\th{phi-gamma}Let $\pi$ be an irreducible representation of $G$, which
appears in the Deligne-Lusztig representation $R_{\theta,w}$
for a quasi-regular character $\theta:T_w\to \qlb^*$. Then
\eq{phi-gamma}
\sum\limits_{g\in G}\tilPhi_\rho(g)\pi(g)=\gam_\rho(\pi)
\end{equation}
\eth

\sssec{}{The basic conjecture}In order to proceed we will
have to assume that the following conjecture holds.

Choose a maximal unipotent subgroup in $\bfU\subset\bfG$.
Denote by $\bfB$ the Borel subgroup of $\bfG$ which normalizes
$\bfU$.
Let $\bfX=\bfG/\bfU$ and let $\bfr:\bfG\to\bfX$ denote the natural
projection. By the definition one has canonical identification
$\bfT\simeq \bfB/\bfU$. Hence $\bfT$ is naturally embedded into
$\bfX$.

\medskip

\noindent
\conj{fini}
The complex $\bfr_! B_\rho$ vanishes outside of $\bfT$.
\econj
One can show that \refco{fini} implies that \refe{phi-gamma}
holds for any $\theta$.

\medskip


\sec{}{Appendix (by V.~Vologodsky)}

\ssec{}{Notations}

In this appendix we assume that $F$ is a finite extension of $\QQ_p$.
We denote by $\grm$ the maximal ideal of $\calO_F$.

Let $\bfC$ be a smooth, quasi-projective curve over $F$,
$c_0$ be a $F$-point of $\bfC$. Let
$\bfC^{\prime}=\bfC\backslash c_0$.

Let $\bff:\bfX\longrightarrow \bfC^{\prime} $
be a smooth, proper morphism of relative dimension $n$ and $\omega$
be a differential form
$\omega\in \Omega ^{\top}_{\bfX/\bfC^{\prime}} $ of the top degree.

For an $F$-point $c\in \bfC^{\prime}(F)$
we denote by $|\omega| $ the measure associated to $\omega$ on the space of
$F$-points of the fiber $\bfX_c$ over $c$. Put $V(c)=\int_{\bfX_c(K)} |\omega| $.
Choose a coordinate $=\bft$ on a neighborhood of $c_0$ with $\bft(c_0)=0$.

\th{}
There exist integers $n$ and $r$ such that $V(tt^{\prime})=V(t)$
for all $t\in \grm^n$,
$t^{\prime}\in (\calO_F^*)^r$.
\eth

This theorem implies easily \reft{vadik}.
\prf
First, we can make use of Nagata's theorem to construct a
proper, integral scheme $\obX$ over $\bfC$ such that
$\obX \times_\bfC \bfC^{\prime}\simeq \bfX$.
By Hironaka's theorem we can blow up $\obX$ to build
$\tbf$: $\bfY\longrightarrow \bfC$
such that the union of the fiber $\bfY_{c_0}$ over $c_0$ and the closure
the zero locus of
the differential form $\omega$ on $\bfY\backslash \bfY_{c_0}$
is a normal crossing
divisor. It follows that for any $F$-point $a$ of $\bfY_{c_0}$    there
exist local coordinates $\bfx_i$ $(i=0,1...n)$ on a neighborhood
$\bfZ_a$
of $a$
satisfying the following conditions:

\medskip
i) $\bft(\tbf) ={\bf u\Pi_i x_i^{n_i}}$,
for some function $\bfu$ on $\bfZ_a$ with $\bfu(a)\ne 0$ and $n_i \geq 0$

\medskip
ii)$\omega = v \Pi_i x_i^{m_i}dx_0dx_1dx_2...dx_{n-1}$ on
$Z_a\backslash Y_{c_0}$ for some function $v$ on $Z_a$ with $v(a)\ne 0$
and $m_i\in \ZZ$.

It is clear that for a point $c\in \bfC^{\prime}(F)$ sufficiently close
to $c_0$ the fibers $\bfX_c(F)$ and $\bfY_c(F)$ have the same measure.
Hence we can replace $\bfX$ by ${\bf Y\backslash Y}_{c_0}$.
By the implicit function theorem we can choose a neighborhood $Z^{\prime}_a$
of $a$ in p-adic topology: $a\in Z^{\prime}_a \subset Z_a$ such that
the coordinates $x_i$ give an isomorphism (of sets)
$Z^{\prime}_a \simeq \grm^l\times \grm^l \times...\times \grm^l$, where
$l$ is a positive integer. Moreover, if we choose $l$ sufficiently large
the function $(\frac{u(a)}{u})^{n_1^{-1}}$ is well defined on $Z_a'$,
hence, changing
$\bfx_1$ for  $(\frac{u(a)}{\bfu})^{n_1^{-1}}\bfx_1$ we can suppose that on
$Z_a^{\prime}$ one has $t(\tilde f) =u(a)\Pi_i x_i^{n_i}$.
Without loss of generality we can assume that $\vert v \vert$ is constant on
$Z_a^{\prime}$.
Since $\bfY_{c_0}(F)$ is compact it can be covered by finitely many open sets
$Z_a^{\prime}$ satisfying the properties stated above.
In fact one can choose them to be disjoint.
To prove it we need the following    lemma.
\lem{}
Let $Q$ be an open, compact subset of the affine space $F^n$.
 We claim that $Q$ can be covered by finitely many disjoint balls
contained in $Q$.
(A ball is a subset of the form $B_{r,a}=\{(x_1,...x_n)\in K^n,
\vert x_i-a_i \vert <r\})$.
\elem
We apply the lemma  to the sets $Z^{\prime}_q\bigcup _{j<q} Z^{\prime}_{a_j}$
$(q=1,...,k)$.
(By construction each of them is identified with an open compact subset
of $F^{n+1}$).

\prf
Since  $Q$ is compact we can cover it by finitely many balls
contained in $Q$. Let $r_0$ be the smallest radius of these balls.
Pick an integer $l$ such that
$-\log_q r_0<l$. Consider the projection $p:Q \to (F/\grm^l)^n$.
The pre-image of each element of Im$p$ is a ball.
These balls constitute the desired covering.

    Now the proposition follows from the following simple lemma

\lem{}
Let $\omega$ stand for differential form
$$
\omega =\prod\limits_i \bfx_i^{m_i}d\bfx_0d\bfx_1d\bfx_2...d\bfx_{n-1}
$$
where $m_i \in \ZZ$ and
 $$
V(t)=\int_{\Pi_i x_i^{n_i}=t;  x_i \in m}
 |\omega|
$$
Then $V(t\cdot t^{\prime})=V(t)$ for any $t\in \grm$, $t^{\prime}
\in (\calO_F^*)^{n_1}$.
\elem
\epr
\epr

\end{document}